\documentclass[10pt]{amsart}

\pdfoutput=1

\usepackage[text={420pt,660pt},centering]{geometry}
\usepackage{amssymb}
\usepackage{bbm}
\usepackage{mathtools}
\usepackage{array}
\usepackage[colorlinks=true,pdfstartview=Fit,linkcolor=blue,citecolor=blue,urlcolor=blue,pagebackref=false]{hyperref}
\usepackage{microtype}
\newtheorem{proposition}{Proposition}
\newtheorem{theorem}[proposition]{Theorem}
\newtheorem{lemma}[proposition]{Lemma}
\newtheorem{corollary}[proposition]{Corollary}
\theoremstyle{definition}
\newtheorem{definition}[proposition]{Definition}
\theoremstyle{remark}
\newtheorem{remark}[proposition]{Remark}
\theoremstyle{plain}

\numberwithin{equation}{section}
\numberwithin{proposition}{section}

\newcommand{\R}{\mathbb{R}}
\newcommand{\C}{\mathbb{C}}
\newcommand{\E}{\mathbb{E}}
\renewcommand{\P}{\mathbb{P}}

\renewcommand{\leq}{\leqslant}
\renewcommand{\geq}{\geqslant}
\renewcommand{\bar}{\overline}
\renewcommand{\d}{\mathrm{d}}
\newcommand{\1}{\mathbf{1}}

\newcommand{\Tr}{\operatorname{Tr}}
\newcommand{\diag}{\operatorname{diag}}

\begin{document}

\title{Free energy of 2-spin spherical vector spin glasses}
\author{Hong-Bin Chen}
\address{NYU--ECNU Institute of Mathematical Sciences, NYU Shanghai, China}
\email{hongbin.chen@nyu.edu}
\author{Yanjin Xiang}
\address{School of Mathematical Sciences, Peking University}
\email{yjxiang@stu.pku.edu.cn}
\author{Zhihua Zhang}
\address{School of Mathematical Sciences, Peking University}
\email{zhzhang@math.pku.edu.cn}
\date{}
\hypersetup{
  pdftitle={Free energy of 2-spin spherical vector spin glasses},
  pdfauthor={Hong-Bin Chen, Yanjin Xiang, Zhihua Zhang}
}

\begin{abstract}
We compute the limit of free energy of 2-spin spherical vector spin glasses with possibly non-convex interaction and couplings satisfying a finite fourth-moment assumption. The couplings need not be Gaussian. The limit is given by a logarithmic-determinant variational formula, equivalently by a finite-dimensional matrix Dyson equation, and its zero-temperature limit yields the ground-state energy. The limiting formulas depend only on the covariance of the couplings. We also study fluctuations: at an interior minimizer, an explicit variance bound transfers the limiting law of a logarithmic linear statistic to the free energy. For Gaussian couplings, this yields Gaussian fluctuations. At a boundary minimizer, assuming a simple regular edge with GOE Tracy--Widom fluctuations and a one-edge reduction, we transfer the edge fluctuations to the free energy.
\end{abstract}

\maketitle
\markboth{Free energy of 2-spin spherical vector spin glasses}{Free energy of 2-spin spherical vector spin glasses}

\section{Introduction}
Throughout, $M_d(\mathbb{F})$ denotes the space of $d\times d$ matrices with entries in $\mathbb{F}$ for some set $\mathbb{F}$; and $\mathbb N:=\{1,2,\ldots\}$. For Hermitian matrices of the same size, or self-adjoint elements of the same $C^*$-algebra, the relation $A\succeq B$ (resp.\ $A\succ B$) means that $A-B$ is positive semidefinite (resp.\ positive definite).

Fix $d\in\mathbb N$. Let $(g_{ij})_{i,j\in\mathbb N}$ be independent copies of a real matrix $g_{11}:=(g_{11}^{ab})_{a,b\in\{1,\ldots,d\}}$ satisfying
\begin{equation}\label{eq:coupling-covariance}
\E g_{11}=0,\qquad
\E g_{11}^{ab}g_{11}^{ce}=\Delta_{ab,ce}\quad(a,b,c,e\in\{1,\ldots,d\}),
\qquad\text{and}\qquad
\E\|g_{11}\|_{\mathrm F}^4<\infty,
\end{equation}
where $\|g\|_{\mathrm F}:=(\sum_{a,b=1}^d|g^{ab}|^2)^\frac{1}{2}$ and $\Delta:=(\Delta_{ab,ce})_{(a,b),(c,e)}\succeq0$. The entries within each $g_{ij}$ need not be independent. We emphasize that we do not assume that $g_{11}$ is Gaussian. The space of vector spherical spins is
\begin{equation*}
\mathcal S_N:=\left\{\sigma=(\sigma^1,\ldots,\sigma^d)\in(\R^N)^d:\|\sigma^a\|^2=N\text{ for every }a\in\{1,\ldots,d\}\right\},
\end{equation*}
and the Hamiltonian is
\begin{equation}\label{eq:hamiltonian}
H_N(\sigma):=\frac1{\sqrt N}\sum_{i,j=1}^N\sum_{a,b=1}^dg_{ij}^{ab}\sigma_i^a\sigma_j^b.
\end{equation}
For $\sigma,\sigma'\in\mathcal S_N$, define the matrix overlap $R(\sigma,\sigma')\in M_d(\R)$ by
\begin{equation*}
R(\sigma,\sigma')_{ac}:=\frac1N\sigma^a\cdot \sigma'^c.
\end{equation*}
It follows from \eqref{eq:coupling-covariance} that
\begin{equation}\label{eq:covariance}
\E H_N(\sigma)H_N(\tau)=N\xi(R(\sigma,\tau)),
\quad\text{where}\quad
\xi(R):=\sum_{a,b,c,e=1}^d\Delta_{ab,ce}R_{ac}R_{be}\quad \text{for $R\in M_d(\R)$}.
\end{equation}
The condition $\Delta\succeq0$ concerns the composite indices $(a,b)$ and $(c,e)$ in \eqref{eq:coupling-covariance}, whereas convexity of $R\mapsto\xi(R)$ concerns the quadratic form in the composite indices $(a,c)$ and $(b,e)$ in \eqref{eq:covariance}. Therefore, $\Delta\succeq0$ does not imply that $R\mapsto\xi(R)$ is convex.

Let $\mu_N$ be the normalized surface measure on the sphere of radius $\sqrt N$ in $\R^N$, and let $\mu_N^{\otimes d}$ be its $d$-fold product. For every $\beta\geq0$, define
\begin{equation*}
F_N(\beta):=\frac1N\log\int_{\mathcal S_N}\exp\bigl(\beta H_N(\sigma)\bigr)\,\d\mu_N^{\otimes d}(\sigma).
\end{equation*}
To state the limit, choose $m\in\mathbb N$ and matrices $A_1,\ldots,A_m\in M_d(\R)$ such that
\begin{equation}\label{eq:factorization}
\Delta_{ab,ce}=\sum_{r=1}^m(A_r)_{ab}(A_r)_{ce}.
\end{equation}
The factorization \eqref{eq:factorization} is only a factorization of the covariance matrix $\Delta$; it imposes no representation assumption on the random matrix $g_{ij}$. Let $(\mathcal A,\tau)$ be a tracial $C^*$-probability space with faithful trace, and write $1$ for the unit of $\mathcal A$. Let $c_1,\ldots,c_m$ be a free circular family, so that
\begin{equation*}
c_r:=\frac{s_{r,1}+\mathrm i s_{r,2}}{\sqrt2},
\end{equation*}
where $(s_{r,k})_{r\in\{1,\ldots,m\},\,k\in\{1,2\}}$ is a freely independent family of variance-one semicircular variables. Thus $\tau(c_r)=0$, $\tau(c_rc_s^*)=\delta_{rs}$, and $\tau(c_rc_s)=0$.

For $B\in M_d(\C)$ and $x\in\mathcal A$, the tensor $B\otimes x\in M_d(\mathcal A)$ is the operator-valued matrix whose $(a,b)$ entry is $B_{ab}x$; when both factors are finite matrices, $\otimes$ denotes their Kronecker product. Define the self-adjoint operator
\begin{equation}\label{eq:free-operator}
\mathsf M:=\frac12\sum_{r=1}^m\left(A_r\otimes c_r+A_r^{\mathsf T}\otimes c_r^*\right) \in M_d(\mathcal A).
\end{equation}
For $X=(X_{ab})_{a,b\in\{1,\ldots,d\}}\in M_d(\mathcal A)$, define the normalized trace
\begin{equation*}
\tau_d(X):=\frac1d\sum_{a=1}^d\tau(X_{aa}).
\end{equation*}
Equivalently, the entries of \eqref{eq:free-operator} and its representation by free semicircular variables are
\begin{align}
\mathsf M_{ab}&=\frac12\sum_{r=1}^m\left((A_r)_{ab}c_r+(A_r)_{ba}c_r^*\right),\notag\\
\mathsf M&=\frac1{2\sqrt2}\sum_{r=1}^m\left((A_r+A_r^{\mathsf T})\otimes s_{r,1}+\mathrm i(A_r-A_r^{\mathsf T})\otimes s_{r,2}\right).\label{eq:free-operator-semicircular}
\end{align}
For $\lambda\in\R^d$, write $D_\lambda:=\diag(\lambda_1,\ldots,\lambda_d)$ and $\1:=(1,\ldots,1)\in\R^d$, and define
\begin{equation}\label{eq:edge-function}
g_\beta(\lambda):=\lambda_{\max}(2\beta\mathsf M-D_\lambda\otimes1).
\end{equation}
Define
\begin{equation}\label{eq:admissible-domain}
\Lambda_\beta
:
=\left\{\lambda\in\R^d:g_\beta(\lambda)<0\right\}.
\end{equation}
The function $g_\beta$ is continuous and satisfies
\begin{equation}\label{eq:edge-shift}
g_\beta(\lambda+t\1)=g_\beta(\lambda)-t
\qquad\text{for every }\lambda\in\R^d\text{ and }t\in\R.
\end{equation}
The equality in \eqref{eq:admissible-domain} and the continuity of $g_\beta$ imply that $\Lambda_\beta$ is open. It is convex because strict positivity of $D_\lambda\otimes1-2\beta\mathsf M$ is preserved under convex combinations. Moreover, \eqref{eq:edge-shift} gives
\begin{equation}\label{eq:admissible-closure-boundary}
\overline{\Lambda}_\beta
=\left\{\lambda\in\R^d:g_\beta(\lambda)\leq0\right\}
\qquad\text{and}\qquad
\partial\Lambda_\beta
=\left\{\lambda\in\R^d:g_\beta(\lambda)=0\right\}.
\end{equation}
For $\lambda\in\Lambda_\beta$, define
\begin{equation}\label{eq:limiting-objective}
\phi_\beta(\lambda)
:=\frac12\sum_{a=1}^d(\lambda_a-1)
-\frac d2\tau_d\log(D_\lambda\otimes1-2\beta\mathsf M).
\end{equation}

\begin{theorem}\label{thm:main}
Under \eqref{eq:coupling-covariance}, for every $\beta>0$, the following limit holds almost surely:
\begin{equation}\label{eq:main-formula}
\lim_{N\to\infty}F_N(\beta)
=\inf_{\lambda\in\Lambda_\beta}\phi_\beta(\lambda).
\end{equation}
If $g_{11}$ is Gaussian, the convergence also holds in $L^p$ for every $1\leq p<\infty$.
\end{theorem}

We write $\operatorname{id}\otimes\tau:M_d(\mathcal A)\longrightarrow M_d(\C)$ for the entrywise application of $\tau$, namely, $(\operatorname{id}\otimes\tau)(X):=(\tau(X_{ab}))_{a,b\in\{1,\ldots,d\}}$ for $X=(X_{ab})\in M_d(\mathcal A)$, where $\operatorname{id}$ denotes the identity map on $M_d(\C)$. Corollary~\ref{cor:finite-dimensional} replaces the operator-valued logarithm in \eqref{eq:main-formula} by a $d\times d$ matrix Dyson equation and gives a finite-dimensional criterion for membership in $\Lambda_\beta$. To state it, define the completely positive covariance map $\eta:M_d(\C)\longrightarrow M_d(\C)$ associated with $\mathsf M$ by
\begin{equation}\label{eq:eta}
\eta(Q):=(\operatorname{id}\otimes\tau)\left[\mathsf M(Q\otimes1)\mathsf M\right]
=\frac14\sum_{r=1}^m\left(A_rQA_r^{\mathsf T}+A_r^{\mathsf T}QA_r\right).
\end{equation}
Equivalently, its entries satisfy $\eta(Q)_{ac}=\frac14\sum_{b,e=1}^d\bigl(\Delta_{ab,ce}+\Delta_{ba,ec}\bigr)Q_{be}$. Consequently, $\eta$ is determined by $\Delta$. By the operator-valued distribution of $\mathsf M$ we mean the collection of all $M_d(\C)$-valued moments $(\operatorname{id}\otimes\tau)[B_0\mathsf M B_1\mathsf M\cdots\mathsf M B_k]$, where $B_0,\ldots,B_k\in M_d(\C)$ are identified with $B_j\otimes1$. Since $\mathsf M$ is an $M_d(\C)$-valued semicircular element, these moments are determined by its covariance map $\eta$. Hence $\Lambda_\beta$, $\phi_\beta$, and the right-hand sides of \eqref{eq:main-formula}, \eqref{eq:finite-dimensional-formula}, and \eqref{eq:ground-state-limit} are independent of the factorization in \eqref{eq:factorization}.

For $\beta>0$ and $\lambda\in\Lambda_\beta$, define $Q_\lambda\in M_d(\C)$ by
\begin{equation}\label{eq:Q-lambda}
Q_\lambda:=(\operatorname{id}\otimes\tau)
\left[(D_\lambda\otimes1-2\beta\mathsf M)^{-1}\right].
\end{equation}
Proposition~\ref{prop:dyson-evaluation} shows that $Q_\lambda\succ0$ and satisfies the matrix Dyson equation
\begin{equation}\label{eq:matrix-dyson}
Q_\lambda=\left(D_\lambda-4\beta^2\eta(Q_\lambda)\right)^{-1}.
\end{equation}
The matrix Dyson equation \eqref{eq:matrix-dyson} can have more than one positive solution. The matrix $Q_\lambda$ defined by \eqref{eq:Q-lambda} is the resolvent solution; the proof of Proposition~\ref{prop:dyson-evaluation} obtains it along the continuous resolvent path beginning at $D_\lambda^{-1}$.

\begin{corollary}\label{cor:finite-dimensional}
For every $\beta>0$, we have
\begin{equation}\label{eq:finite-dimensional-formula}
\lim_{N\to\infty}F_N(\beta)
=\inf_{\lambda\in\Lambda_\beta}
\left\{
\frac12\sum_{a=1}^d(\lambda_a-1)
+\frac12\log\det Q_\lambda
-\beta^2\Tr\bigl(Q_\lambda\eta(Q_\lambda)\bigr)
\right\}.
\end{equation}
Moreover, for every $\lambda\in\R^d$, we have
\begin{equation}\label{eq:admissible-finite}
\lambda\in\Lambda_\beta
\quad\Longleftrightarrow\quad
\max_{\substack{S\in M_d(\C),\ S=S^*,\ S\succeq0\\\Tr S=1}}
\left\{-\Tr(D_\lambda S)+4\beta\Tr\left[(S^\frac{1}{2}\eta(S)S^\frac{1}{2})^\frac{1}{2}\right]\right\}<0.
\end{equation}
\end{corollary}

The infimum in \eqref{eq:main-formula} need not be attained in $\Lambda_\beta$. Correspondingly, the finite-dimensional minimizing precision matrices need not have a spectral gap bounded uniformly away from zero.

Theorem~\ref{thm:ground-state} gives the zero-temperature limit of Theorem~\ref{thm:main}. For every $\beta>0$, write $F(\beta):=\lim_{N\to\infty}F_N(\beta)$.

\begin{theorem}\label{thm:ground-state}
Under \eqref{eq:coupling-covariance}, almost surely,
\begin{equation}\label{eq:ground-state-limit}
\lim_{\beta\to\infty}\frac{F(\beta)}\beta
=\lim_{N\to\infty}\frac1N\max_{\sigma\in\mathcal S_N}H_N(\sigma)
=2\max_{\substack{Q\in M_d(\R),\ Q=Q^{\mathsf T},\ Q\succeq0\\Q_{aa}=1\text{ for every }a\in\{1,\ldots,d\}}}
\Tr\left[(Q^\frac{1}{2}\eta(Q)Q^\frac{1}{2})^\frac{1}{2}\right].
\end{equation}
If $g_{11}$ is Gaussian, the convergence of $N^{-1}\max_{\sigma\in\mathcal S_N}H_N(\sigma)$ in \eqref{eq:ground-state-limit} also holds in $L^p$ for every $1\leq p<\infty$.
\end{theorem}

The limits in Theorems~\ref{thm:main} and~\ref{thm:ground-state} and Corollary~\ref{cor:finite-dimensional} depend on the law of $g_{11}$ only through $\Delta$. The assumption $\E\|g_{11}\|_{\mathrm F}^4<\infty$ does not imply the $L^p$ conclusions for every finite $p$.

Universality for spherical spins requires control of configurations with large coordinates: here $\sup_{\sigma\in\mathcal S_N}\max_{a,i}|\sigma_i^a|=\sqrt N$, so an argument requiring an $N$-independent coordinate bound does not apply. For scalar spherical mixed $p$-spin models, Sawhney and Sellke \cite[Theorem 1.5]{sawhney-sellke} prove almost-sure universality of the limiting free energy and ground-state energy under a finite $2p$-th moment assumption on the $p$-spin couplings. In the present 2-spin vector model, universality follows from convergence of the spectra and normalized traces of the associated block matrices, established in Proposition~\ref{prop:strong-convergence} by the matrix comparison and truncation arguments of Appendix~\ref{app:fourth-moment}. These conclusions concern the limiting free energy and ground-state energy; fluctuation limits require the additional hypotheses stated below.

\smallskip

\noindent \textbf{Fluctuations.}
For $\lambda\in\partial\Lambda_\beta$, the operator $D_\lambda\otimes1-2\beta\mathsf M$ is positive semidefinite but not strictly positive, so \eqref{eq:limiting-objective} need not be defined. For $\lambda\in\R^d$ such that $\log|D_\lambda\otimes1-2\beta\mathsf M|$ is $\tau_d$-integrable,\footnote{For a self-adjoint $B\in M_d(\mathcal A)$, let $\nu_B$ be its spectral distribution under $\tau_d$, characterized by $\tau_d(f(B))=\int_{\R}f(x)\,\d\nu_B(x)$ for continuous $f$. Here integrability means $\int_{\R}|\log|x||\,\d\nu_B(x)<\infty$, and $\tau_d(\log|B|):=\int_{\R}\log|x|\,\d\nu_B(x)$.} define
\begin{equation}\label{eq:two-sided-objective}
\widehat\phi_\beta(\lambda)
:=\frac12\sum_{a=1}^d(\lambda_a-1)
-\frac d2\tau_d\!\left[\log\left|D_\lambda\otimes1-2\beta\mathsf M\right|\right],
\end{equation}
where $|B|:=(B^*B)^\frac{1}{2}=(B^2)^\frac{1}{2}$ for self-adjoint $B\in M_d(\mathcal A)$, with the positive square root defined by functional calculus. For every $\lambda\in\Lambda_\beta$, we have $\widehat\phi_\beta(\lambda)=\phi_\beta(\lambda)$. When the spectral density vanishes at square-root rate at an edge, the logarithmic singularity in \eqref{eq:two-sided-objective} is integrable.

\begin{definition}[Boundary multiplier]\label{def:boundary-multiplier}
Fix $\beta>0$ and $\lambda_\beta\in\partial\Lambda_\beta$. Assume that $\widehat\phi_\beta$ is continuously differentiable in a neighborhood of $\lambda_\beta$, that
\begin{equation*}
\widehat\phi_\beta(\lambda_\beta)=F(\beta),
\end{equation*}
and that $g_\beta$ is differentiable at $\lambda_\beta$. The unique number $\theta_\beta\geq0$ satisfying
\begin{equation}\label{eq:boundary-multiplier}
\nabla\widehat\phi_\beta(\lambda_\beta)
=-\theta_\beta\nabla g_\beta(\lambda_\beta),
\end{equation}
is called the \textit{boundary multiplier} associated with $\lambda_\beta$.
\end{definition}
Recall $Q_\lambda$ from \eqref{eq:Q-lambda}. Equation \eqref{eq:edge-shift} gives $\lambda_\beta+t\1\in\Lambda_\beta$ for every $t>0$. Write
\begin{equation}\label{eq:boundary-resolvent}
Q^\partial_{\lambda_\beta}:=\lim_{t\downarrow0}Q_{\lambda_\beta+t\1},
\end{equation}
whenever this limit exists in matrix norm.

\begin{proposition}\label{prop:boundary-resolvent-multiplier}
Under the hypotheses of Definition~\ref{def:boundary-multiplier}, the multiplier $\theta_\beta$ is well defined, $Q^\partial_{\lambda_\beta}\in M_d(\C)$ exists, and
\begin{equation}\label{eq:boundary-multiplier-trace}
\theta_\beta
=\nabla\widehat\phi_\beta(\lambda_\beta)\cdot\1
=\frac12\left(d-\Tr Q^\partial_{\lambda_\beta}\right)\geq0.
\end{equation}
In particular, $\theta_\beta>0$ if and only if $\Tr Q^\partial_{\lambda_\beta}<d$.
\end{proposition}
The proof is given at the start of Subsection~\ref{sec:boundary-fluctuations}.

At an interior minimizer, Theorem~\ref{thm:interior-fluctuations} transfers any limiting law of the centered logarithmic linear statistic to $N(F_N(\beta)-F(\beta))$ after deterministic centering, under \eqref{eq:interior-resolvent-variance}. The table summarizes the resulting Gaussian and boundary conclusions. The labels ``high'', ``low'', and ``critical temperature'' refer to conditions on the variational minimizer at a fixed $\beta$, as specified in the second column. They reflect the usual interpretation of small and large $\beta$ as high and low temperature, respectively; they do not assert a partition of $(0,\infty)$ by a single critical value of $\beta$.
\begin{table}[t]
\centering
\small
\setlength{\extrarowheight}{4pt}
\begin{tabular}{@{}>{\raggedright\arraybackslash}p{92pt}>{\raggedright\arraybackslash}p{158pt}>{\raggedright\arraybackslash}p{146pt}@{}}
\hline
\textbf{Case} & \textbf{Conditions} & \textbf{Conclusion} \\[4pt]
\hline
``High temperature'' & There is a unique $\lambda_\beta\in\Lambda_\beta$ with $\phi_\beta(\lambda_\beta)=F(\beta)$. & If $g_{11}$ is Gaussian, the free energy has Gaussian fluctuations at scale $N^{-1}$ after deterministic centering (Corollary~\ref{cor:gaussian-interior-fluctuations}). \\[4pt]
\hline
``Low temperature'' & $\lambda_\beta\in\partial\Lambda_\beta$, $\widehat\phi_\beta(\lambda_\beta)=F(\beta)$, and $\theta_\beta>0$; $(g_\beta,\widehat\phi_\beta)$ is nondegenerate at $\lambda_\beta$ (Definition~\ref{def:boundary-nondegeneracy}), $\lambda_\beta$ exposes a simple regular GOE edge (Definition~\ref{def:simple-regular-goe-edge}), and \eqref{eq:single-edge-reduction} holds. & The assumed edge GOE Tracy--Widom law is transferred to the centered free energy at scale $N^{-2/3}$ (Theorem~\ref{thm:simple-regular-goe-free-energy}). \\[4pt]
\hline
``Critical temperature'' & $\lambda_\beta\in\partial\Lambda_\beta$, $\widehat\phi_\beta(\lambda_\beta)=F(\beta)$, and $\theta_\beta=0$. & No fluctuation limit is proved here. \\[4pt]
\hline
Other boundary cases & For a boundary minimizer $\lambda_\beta$, $(g_\beta,\widehat\phi_\beta)$ is not nondegenerate at $\lambda_\beta$ in the sense of Definition~\ref{def:boundary-nondegeneracy}, a condition in Definition~\ref{def:simple-regular-goe-edge} fails, or \eqref{eq:single-edge-reduction} fails. & No fluctuation limit is proved here. Failure of these hypotheses does not by itself identify a critical temperature. \\[4pt]
\hline
\end{tabular}
\end{table}

For all sufficiently small $\beta>0$, there exists a unique $\lambda_\beta\in\Lambda_\beta$ such that $\phi_\beta(\lambda_\beta)=F(\beta)$. The implicit function theorem at $\beta=0$, where $\phi_0$ has the nondegenerate minimizer $\1$ with $\nabla^2\phi_0(\1)=I_d/2$, gives a critical point in $\Lambda_\beta$ for all sufficiently small $\beta>0$; Lemma~\ref{lem:interior-minimizer-uniqueness} shows that this point is the unique minimizer.

Under only $d\geq1$ and $\Delta\succeq0$, we do not prove that there exists $\beta_c>0$ such that
\begin{equation*}
\left\{\beta>0:\exists\,\lambda_\beta\in\Lambda_\beta,\ \phi_\beta(\lambda_\beta)=F(\beta)\right\}
=(0,\beta_c),
\end{equation*}
and such that, for every $\beta>\beta_c$, there exists $\lambda_\beta\in\partial\Lambda_\beta$ with $\widehat\phi_\beta(\lambda_\beta)=F(\beta)$ and $\theta_\beta>0$. Lemma~\ref{lem:interior-minimizer-uniqueness} shows that there is at most one $\lambda_\beta\in\overline\Lambda_\beta$ with finite value $\widehat\phi_\beta(\lambda_\beta)=F(\beta)$. Theorem~\ref{thm:main} and Section~\ref{sec:fluctuations} do not exclude more than one eigenvalue of $2\beta M_N-D_{\lambda_\beta}\otimes I_N$ at distance $O_{\mathbb P}(N^{-2/3})$ from $0$, or more than one value of $\beta$ at which the condition $\lambda_\beta\in\Lambda_\beta$ changes to $\lambda_\beta\in\partial\Lambda_\beta$. Thus the table applies separately at each fixed $\beta$.

For the equal-size bipartite spherical model with Gaussian couplings in \eqref{eq:bipartite-specialization}, \cite[Theorems 2.2(i) and 2.4(i)]{baik-lee} treats $0<\beta<1$ with $\lambda_\beta\in\Lambda_\beta$, while Section~\ref{sec:bipartite-fluctuations} verifies the hypotheses of Theorem~\ref{thm:simple-regular-goe-free-energy} for $\beta>1$, with $\lambda_\beta=(2\beta,2\beta)\in\partial\Lambda_\beta$ and $\theta_\beta=(\beta-1)/\beta>0$. The value $\beta=1$ is excluded from both noncritical regimes.

At zero temperature, Theorem~\ref{thm:ground-state-fluctuations} gives the $N^{-2/3}$-scale GOE Tracy--Widom limit in \eqref{eq:ground-state-tw}. A universal single Tracy--Widom law cannot hold for every $\Delta\succeq0$: if the Hamiltonian is a sum of independent scalar spherical Hamiltonians, then the low-temperature free energy and the ground-state energy contain sums of independent edge contributions.

\smallskip

\noindent \textbf{Related works.}
The scalar spherical Sherrington--Kirkpatrick model was introduced by Kosterlitz, Thouless, and Jones \cite{kosterlitz-thouless-jones}. The Crisanti--Sommers representation originated in \cite{crisanti-sommers}. Talagrand proved the corresponding variational formula for mixed even interactions in \cite[Theorem 1.1]{talagrand-spherical}, and Chen extended the formula to general mixed interactions in \cite[Theorem 1]{chen-spherical}. In the pure 2-spin model without an external field, diagonalization reduces the free energy to the spectrum of one random matrix. Baik and Lee identified the Gaussian and GOE Tracy--Widom fluctuation limits on the two sides of the critical inverse temperature in \cite[Theorem 1.1]{baik-lee-spherical}. Landon analyzed the critical inverse temperature and its transition window in \cite[Theorems 1.1 and 1.2]{landon-critical}; Johnstone, Klochkov, Onatski, and Pavlyshyn identified the limiting law throughout the fixed-parameter critical window in \cite[Theorem 1]{johnstone-klochkov-onatski-pavlyshyn}.

For vector spins and coupled spherical systems, Panchenko and Talagrand initiated the study of overlap-constrained multiple spherical models in \cite{panchenko-talagrand}. Panchenko proved a Parisi formula for mixed even interactions with vector spins in \cite[Theorem 1]{panchenko-vector}. Ko identified the constrained spherical free energy in \cite[Theorem 2.1]{ko} and derived a vector Crisanti--Sommers representation in \cite[Theorems 1 and 2]{ko-crisanti-sommers}. For the spherical mixed vector $p$-spin model, Auffinger and Zhou studied the Parisi functional and its zero-temperature limit in \cite[Theorem 1 and Proposition 4]{auffinger-zhou}. For 2-spin shared disorder with diagonal component weights, the fixed-overlap ground-state formulas in \cite[Proposition 4]{auffinger-zhou} and \cite[Equation (1.7) and Theorem 1.2]{belius-frober-ko} agree with the restriction of \eqref{eq:ground-state-limit}, because $\xi'(Q)=4\eta(Q)$ in that subclass.

For multi-species spherical models, Bates and Sohn proved the Parisi formula under the convexity assumption on the overlap covariance in \cite[Theorem 1.3]{bates-sohn} and obtained its Crisanti--Sommers representation in \cite[Theorem 2.2 and Corollary 2.3]{bates-sohn-crisanti-sommers}. Subag developed the corresponding TAP representation and identified the Onsager correction in \cite[Theorems 2 and 3]{subag-tap-i}. Conditional on the convergence hypothesis in \cite[Equation (1.6)]{subag-tap-ii}, Subag evaluated the pure multi-species free energy for $p\geq3$ in \cite[Theorem 2]{subag-tap-ii} and treated $p=2$ in \cite[Appendix]{subag-tap-ii}. For balanced models, the non-convex lower bound of \cite[Theorem 1.3]{bates-sohn-balanced} is an equality in several subclasses; the pure bipartite spherical ground-state case is given by \cite[Corollary 2.7]{bates-sohn-balanced}.

The closest non-convex 2-spin example is the bipartite spherical model. Auffinger and Chen obtained a high-temperature variational formula for mixed bipartite interactions in \cite[Theorem 1]{auffinger-chen-bipartite}. For the pure bipartite model, Baik and Lee computed the limiting free energy and its fluctuations on both sides of the critical inverse temperature in \cite[Theorems 2.1 and 2.2]{baik-lee}; Collins-Woodfin and Le identified the critical-window fluctuations in \cite[Theorem 1.1]{collins-woodfin-le}. Write $H^{\mathrm{bip}}_{N,N}$ and $F^{\mathrm{bip}}_{N,N}$ for the equal-size Hamiltonian and free energy in \cite[Equations (1.3) and (2.4)]{baik-lee}. If $g_{11}$ is Gaussian, $d=2$, $\Delta_{12,12}=1$, and every other entry of $\Delta$ is zero, then
\begin{equation}\label{eq:bipartite-specialization}
H_N=\sqrt2H^{\mathrm{bip}}_{N,N}
\qquad\text{and}\qquad
F_N(\beta)=2F^{\mathrm{bip}}_{N,N}(\sqrt2\beta).
\end{equation}
The assumptions for the limiting formulas are $d\in\mathbb N$ and \eqref{eq:coupling-covariance}; no convexity assumption is imposed on $R\mapsto\xi(R)$ in \eqref{eq:covariance}, and the entries of $\Delta$ may correlate different coupling blocks. General vector interactions with non-convex overlap covariance were related to a Hamilton--Jacobi critical-value problem in \cite[Theorem 1.2]{chen-mourrat}, while the centered Ising multi-species free energy was identified without a convexity assumption in \cite[Theorem 1.1]{chen-issa-mourrat}.

Corollary~\ref{cor:gaussian-interior-fluctuations} uses the real second-order Wick expansion and freeness of independent Ginibre matrices and their transposes from \cite[Sections 4--6]{redelmeier}, together with the covariance and central-limit results in \cite[Theorem 6(b), Theorem 11, and Proposition 12]{diaz-mingo}. The block formula in \cite[Theorem 2]{diaz-mingo-belinschi} gives the complex Gaussian second-order Cauchy transform. For Definition~\ref{def:simple-regular-goe-edge}, the regular-edge notion and GOE edge universality follow from \cite[Definition 2.4 and Theorem 2.11]{alt-erdos-kruger-schroder} under the hypotheses of the latter theorem, while \cite[Theorem 2.6]{erdos-henheik-riabov} gives the required block-profile estimate under its hypotheses. Structured models with deterministic zero blocks, including the bipartite linearization, can instead be verified directly from their singular-value representation.

For the first-order formulas, Proposition~\ref{prop:strong-convergence} proves convergence of the spectra and normalized traces of the matrices in \eqref{eq:affine-pencils}. Its proof in Appendix~\ref{app:fourth-moment} uses the comparison inequalities of Brailovskaya and van Handel \cite[Theorems 2.6 and 2.9 and Lemma 9.20]{brailovskaya-van-handel}, applied to independent matrix summands with the same entry covariances, and removes truncation using Bai and Yin \cite[Theorem 2.1, with $k=1$]{bai-yin-products}. Xiang and Zhang \cite[Theorem 3.3 and Remark 3.6]{xiang-zhang} identify the Gaussian comparison limit. This proves the required convergence even when the coordinates within $g_{ij}$ are dependent. We use Kunisky \cite[Theorem 1.9]{kunisky} for the finite-dimensional form of the semicircular spectral edge. Starting from Lehner \cite[Corollary 1.5]{lehner}, one can also obtain the reduced dual formula by introducing a density matrix for the largest eigenvalue and applying Sion's minimax theorem \cite[Theorem 4.2]{sion}; see the discussion following \cite[Theorem 1.9]{kunisky}.

\smallskip

\noindent \textbf{Organization.}
Section~\ref{sec:limits} proves the limiting free-energy and ground-state formulas. Its subsections establish the spectral variational principle, the block random-matrix limit, the matrix Dyson equation, the determinant reduction, and the zero-temperature evaluation. Section~\ref{sec:fluctuations} first proves the conditional fluctuation transfer under $\lambda_\beta\in\Lambda_\beta$ and its Gaussian specialization, then proves the one-edge transfer theorem under $\lambda_\beta\in\partial\Lambda_\beta$ and $\theta_\beta>0$, then proves the zero-temperature regular-edge theorem, and finally verifies the two positive-temperature cases in the Gaussian bipartite model of Baik and Lee. Appendix~\ref{app:fourth-moment} proves Proposition~\ref{prop:strong-convergence}.

\smallskip

\noindent \textbf{Acknowledgements}.
HBC acknowledges funding from the NYU Shanghai Start-Up Fund and support from the NYU--ECNU Institute of Mathematical Sciences at NYU Shanghai. The authors prepared this paper with assistance from large language models, including for suggesting arguments and contributing to drafting, revision, and mathematical checks.

\section{Limiting free energy and ground-state energy}\label{sec:limits}

\subsection{A spectral variational principle}

For $a\in\{1,\ldots,d\}$, let $E_{aa}\in M_d(\R)$ be the matrix whose $(a,a)$ entry is $1$ and whose other entries are zero, and set $P_a:=E_{aa}\otimes I_N$. Identify $\mathcal S_N$ with vectors $x\in\R^{dN}$ satisfying $\|P_ax\|^2=N$ for every $a\in\{1,\ldots,d\}$, where we write $I_N$ for the $N\times N$ identity matrix and recall that $1$ is the unit of $\mathcal A$. We use $\|\cdot\|$ for the Euclidean norm of vectors and the operator norm of matrices or operators.

Let $M_N\in M_{dN}(\R)$ be symmetric and let $\mathsf M\in M_d(\mathcal A)$ be self-adjoint, where $\tau$ is faithful. For $B=B^*\in M_d(\C)$, $R\in M_d(\C)$, and $t\in\R$, define
\begin{equation}\label{eq:affine-pencils}
\begin{aligned}
L_N(B,R,t)&:=B\otimes I_N+t(R\otimes I_N)M_N(R^*\otimes I_N),\\
L(B,R,t)&:=B\otimes1+t(R\otimes1)\mathsf M(R^*\otimes1).
\end{aligned}
\end{equation}
For nonempty compact sets $S,T\subset\R$, write
\begin{equation*}
d_{\mathrm H}(S,T):=\max\left\{\sup_{x\in S}\inf_{y\in T}|x-y|,\sup_{y\in T}\inf_{x\in S}|x-y|\right\}.
\end{equation*}
We now record the spectral variational principle used in Lemma~\ref{lem:rank-one-extraction}.

\begin{lemma}\label{lem:spectral-dual}
Suppose that, for every $B,R,t$ as above and every continuous $f:\R\to\C$,
\begin{gather}
\lim_{N\to\infty}d_{\mathrm H}\bigl(\operatorname{spec}L_N(B,R,t),\operatorname{spec}L(B,R,t)\bigr)=0,\label{eq:affine-spectrum}\\
\lim_{N\to\infty}\frac1{dN}\Tr f(L_N(B,R,t))=\tau_d(f(L(B,R,t))).\label{eq:affine-trace}
\end{gather}
Then, we have
\begin{align}\label{eq:spectral-dual}
\lim_{N\to\infty}\max\left\{u^{\mathsf T}M_Nu:\|P_au\|^2=1,\,\forall a\in\{1,\ldots,d\}\right\}
=\inf_{\lambda\in\R^d}\left\{\sum_{a=1}^d\lambda_a+d\lambda_{\max}(\mathsf M-D_\lambda\otimes1)\right\}.
\end{align}
\end{lemma}

\begin{proof}
\smallskip
\noindent\emph{Step 1: Upper bound in \eqref{eq:spectral-dual}.}
By taking $(B,R,t)=(0,I_d,1)$ in \eqref{eq:affine-spectrum}, we get $\sup_N\|M_N\|<\infty$. On each compact set of coefficients, the matrices in \eqref{eq:affine-pencils} are therefore uniformly bounded and uniformly Lipschitz in $(B,R,t)$. Weyl's inequality and uniform continuity of $f$ on a common spectral interval imply that \eqref{eq:affine-spectrum} and \eqref{eq:affine-trace} are locally uniform in $(B,R,t)$, by a finite-net argument. In particular, these limits hold along convergent coefficient sequences. For every admissible $u$ and every $\lambda\in\R^d$, we have
\begin{equation*}
u^{\mathsf T}M_Nu
\leq\sum_{a=1}^d\lambda_a
+d\lambda_{\max}(M_N-D_\lambda\otimes I_N).
\end{equation*}
Equation \eqref{eq:affine-spectrum}, with $(B,R,t)=(-D_\lambda,I_d,1)$, gives
\begin{equation*}
\lim_{N\to\infty}\lambda_{\max}(M_N-D_\lambda\otimes I_N)
=\lambda_{\max}(\mathsf M-D_\lambda\otimes1).
\end{equation*}
This proves the upper bound in \eqref{eq:spectral-dual}.

\smallskip
\noindent\emph{Step 2: Finite-$t$ balancing and Gaussian rounding.}
For the reverse bound, fix $t>0$ and define
\begin{equation*}
\psi_{N,t}(\lambda):=\sum_{a=1}^d\lambda_a+
\frac d t\log\left[\frac1{dN}\Tr e^{t(M_N-D_\lambda\otimes I_N)}\right].
\end{equation*}
This function is invariant under $\lambda\mapsto\lambda+c\1$. Minimize it on $\sum_a\lambda_a=0$. Set $C_M:=\sup_N\|M_N\|<\infty$. Jensen's inequality, applied to the $N$ coordinate vectors in a block on which $\min_a\lambda_a$ is attained, gives
\begin{equation}\label{eq:gibbs-coercive}
\psi_{N,t}(\lambda)\geq-d(C_M+\min_a\lambda_a)-\frac d t\log d,
\qquad\text{and}\qquad
\psi_{N,t}(0)\leq dC_M.
\end{equation}
Consequently, we have
\begin{equation}\label{eq:gibbs-sublevel}
\psi_{N,t}(\lambda)\leq \psi_{N,t}(0)
\quad\Longrightarrow\quad
-\min_a\lambda_a\leq2C_M+\frac{\log d}{t},
\qquad\text{and}\qquad
\max_a\lambda_a\leq-(d-1)\min_a\lambda_a.
\end{equation}
The zero-sum sublevel in \eqref{eq:gibbs-sublevel} is closed and bounded. Hence a minimizer $\lambda_{N,t}$ exists. For every fixed $t$, all these minimizers lie in one compact set.

Set
\begin{equation*}
J_{N,t}:=M_N-D_{\lambda_{N,t}}\otimes I_N,
\qquad\text{and}\qquad
\rho_{N,t}:=\frac{e^{tJ_{N,t}}}{\Tr e^{tJ_{N,t}}}.
\end{equation*}
Duhamel's formula and cyclicity of the trace give
\begin{equation*}
\frac{\partial \psi_{N,t}}{\partial\lambda_a}(\lambda_{N,t})
=1-d\Tr(P_a\rho_{N,t}).
\end{equation*}
At a minimizer on $\sum_a\lambda_a=0$, the gradient is parallel to $\1$. Its components sum to zero, so
\begin{equation}\label{eq:balanced-rho}
\Tr(P_a\rho_{N,t})=\frac1d
\quad\text{for every }a\in\{1,\ldots,d\}.
\end{equation}
Equation \eqref{eq:gibbs-sublevel} implies $\sup_N\|J_{N,t}\|<\infty$. Since $J_{N,t}$ is real symmetric, so are $\rho_{N,t}$ and $\rho_{N,t}^\frac{1}{2}$. Moreover, we have
\begin{equation}\label{eq:rho-small}
\|\rho_{N,t}\|\leq\frac{e^{2t\|J_{N,t}\|}}{dN}=O_t(N^{-1}).
\end{equation}

Let $y_{N,t}$ be a centered real Gaussian vector with covariance $\rho_{N,t}$. For every real symmetric $C$, Wick's formula gives
\begin{equation}\label{eq:rounding-variance}
\begin{aligned}
\E_w y_{N,t}^{\mathsf T}Cy_{N,t}&=\Tr(C\rho_{N,t}),\\
\operatorname{Var}_w(y_{N,t}^{\mathsf T}Cy_{N,t})
&=2\|\rho_{N,t}^\frac{1}{2}C\rho_{N,t}^\frac{1}{2}\|_{\mathrm{HS}}^2
\stackrel{\text{Hilbert--Schmidt}}{\leq}2\|\rho_{N,t}\|\Tr(C^2\rho_{N,t})
\leq2\|C\|^2\|\rho_{N,t}\|.
\end{aligned}
\end{equation}
Apply \eqref{eq:rounding-variance} to $P_1,\ldots,P_d$ and $J_{N,t}$. By \eqref{eq:balanced-rho}, the first $d$ means equal $d^{-1}$. Define
\begin{equation*}
\varepsilon_{N,t}:=\max\left\{
\max_{a\in\{1,\ldots,d\}}\operatorname{Var}_w(\|P_ay_{N,t}\|^2),
\operatorname{Var}_w(y_{N,t}^{\mathsf T}J_{N,t}y_{N,t})
\right\}.
\end{equation*}
By \eqref{eq:rho-small},
\begin{equation*}
\lim_{N\to\infty}\varepsilon_{N,t}=0.
\end{equation*}
For $\varepsilon_{N,t}>0$, Chebyshev's inequality gives
\begin{align*}
\P_w\left(
\max_{a\in\{1,\ldots,d\}}\left|\|P_ay_{N,t}\|^2-d^{-1}\right|>\varepsilon_{N,t}^{1/4}
\ \text{or}\ 
\left|y_{N,t}^{\mathsf T}J_{N,t}y_{N,t}-\Tr(J_{N,t}\rho_{N,t})\right|>\varepsilon_{N,t}^{1/4}
\right)
\leq(d+1)\varepsilon_{N,t}^\frac{1}{2}.
\end{align*}
If $\varepsilon_{N,t}=0$, all these random variables equal their means almost surely. We may therefore choose a deterministic realization such that
\begin{equation}\label{eq:rounding-realization}
\|P_ay_{N,t}\|^2=\frac1d+o(1),
\qquad\text{and}\qquad
y_{N,t}^{\mathsf T}J_{N,t}y_{N,t}=\Tr(J_{N,t}\rho_{N,t})+o(1).
\end{equation}
After rescaling the blocks, set
\begin{equation*}
u_{N,t}:=\sum_{a=1}^d\frac{P_ay_{N,t}}{\|P_ay_{N,t}\|}.
\end{equation*}
Then $\|P_au_{N,t}\|=1$ and
\begin{equation*}
\|u_{N,t}-\sqrt d\,y_{N,t}\|^2
=\sum_{a=1}^d\left(1-\sqrt d\,\|P_ay_{N,t}\|\right)^2=o(1).
\end{equation*}
Since $\sup_N\|J_{N,t}\|<\infty$, \eqref{eq:rounding-realization} gives
\begin{equation}\label{eq:rounded-spectral-energy}
u_{N,t}^{\mathsf T}J_{N,t}u_{N,t}=d\Tr(J_{N,t}\rho_{N,t})+o(1).
\end{equation}
The entropy identity for $\rho_{N,t}$ implies
\begin{equation}\label{eq:gibbs-entropy}
\frac1t\log\left[\frac1{dN}\Tr e^{tJ_{N,t}}\right]
=\Tr(J_{N,t}\rho_{N,t})+
\frac{-\Tr(\rho_{N,t}\log\rho_{N,t})-\log(dN)}t
\leq\Tr(J_{N,t}\rho_{N,t}).
\end{equation}
Since
\begin{equation*}
u_{N,t}^{\mathsf T}M_Nu_{N,t}
=\sum_{a=1}^d(\lambda_{N,t})_a+u_{N,t}^{\mathsf T}J_{N,t}u_{N,t},
\end{equation*}
equations \eqref{eq:rounded-spectral-energy} and \eqref{eq:gibbs-entropy} give
\begin{equation}\label{eq:spectral-lower-finite-t}
\max_{\|P_au\|=1}u^{\mathsf T}M_Nu
\geq\inf_{\lambda\in\R^d}\psi_{N,t}(\lambda)-o(1).
\end{equation}

\smallskip
\noindent\emph{Step 3: The limits $N\to\infty$ and $t\to\infty$.}
For fixed $t$ and $\lambda$, define
\begin{equation*}
\psi_t(\lambda):=\sum_{a=1}^d\lambda_a+
\frac d t\log\tau_d e^{t(\mathsf M-D_\lambda\otimes1)}.
\end{equation*}
Equation \eqref{eq:affine-trace}, applied to the exponential, gives
\begin{equation*}
\lim_{N\to\infty}\psi_{N,t}(\lambda)=\psi_t(\lambda).
\end{equation*}
Duhamel's formula gives
\begin{equation}\label{eq:gibbs-lipschitz}
\max\left\{
\left|\psi_{N,t}(\lambda)-\psi_{N,t}(\mu)\right|,
\left|\psi_t(\lambda)-\psi_t(\mu)\right|
\right\}
\leq d\|\lambda-\mu\|_1.
\end{equation}
The analogue of \eqref{eq:gibbs-coercive}, obtained from the positive states $X\mapsto\tau(X_{aa})$, and the bound $\psi_t(0)\leq d\|\mathsf M\|$ place the finite and limiting minimizers in the compact set determined by \eqref{eq:gibbs-sublevel}. Pointwise convergence and \eqref{eq:gibbs-lipschitz}, applied to a finite $\varepsilon$-net of that compact set, therefore give
\begin{equation}\label{eq:gibbs-inf-limit}
\lim_{N\to\infty}\inf_\lambda \psi_{N,t}(\lambda)=\inf_\lambda \psi_t(\lambda).
\end{equation}
Faithfulness of $\tau_d$ gives
\begin{equation}\label{eq:laplace-edge}
\lim_{t\to\infty}\frac1t\log\tau_d(e^{tB})=\lambda_{\max}(B),
\end{equation}
for every self-adjoint $B$. For $t\geq1$, \eqref{eq:gibbs-sublevel} places every minimizer of $\psi_t$ in one compact subset of $\sum_a\lambda_a=0$. Both
\begin{equation*}
\lambda\longmapsto\frac1t\log\tau_d e^{t(\mathsf M-D_\lambda\otimes1)}
\quad\text{and}\quad
\lambda\longmapsto\lambda_{\max}(\mathsf M-D_\lambda\otimes1)
\end{equation*}
are $\|\cdot\|_\infty$-Lipschitz. A finite-net argument makes \eqref{eq:laplace-edge} uniform on the compact set determined by \eqref{eq:gibbs-sublevel}. Consequently, we have
\begin{equation}\label{eq:gibbs-to-edge}
\lim_{t\to\infty}\inf_\lambda \psi_t(\lambda)
=\inf_\lambda\left\{\sum_{a=1}^d\lambda_a+d\lambda_{\max}(\mathsf M-D_\lambda\otimes1)\right\}.
\end{equation}
For each fixed $t$, we have
\begin{equation*}
\liminf_{N\to\infty}\max_{\|P_au\|=1}u^{\mathsf T}M_Nu
\stackrel{\eqref{eq:spectral-lower-finite-t},\eqref{eq:gibbs-inf-limit}}{\geq}
\inf_\lambda \psi_t(\lambda).
\end{equation*}
Letting $t\to\infty$ and using \eqref{eq:gibbs-to-edge} proves the lower bound without interchanging the two limits.
\end{proof}

\subsection{The block random-matrix limit}

Let $E_{ij}\in M_N(\R)$ be the matrix unit whose $(i,j)$ entry is $1$ and whose other entries are zero. By \eqref{eq:hamiltonian}, if $x:=(\sigma^1,\ldots,\sigma^d)$, then
\begin{equation}\label{eq:block-pencil}
M_N:=\frac1{2\sqrt N}\sum_{i,j=1}^N\bigl(g_{ij}\otimes E_{ij}+g_{ij}^{\mathsf T}\otimes E_{ji}\bigr),
\qquad\text{and}\qquad
H_N(\sigma)=x^{\mathsf T}M_Nx.
\end{equation}
In particular, $(M_N)_{(a,i),(b,j)}=(g_{ij}^{ab}+g_{ji}^{ba})/(2\sqrt N)$. In the scalar unit-variance case, the off-diagonal variance is $1/(2N)$ and the limiting spectral edge is $\sqrt2$.

When $g_{11}$ is Gaussian, \eqref{eq:factorization} permits the realization $g_{ij}=\sum_{r=1}^m A_rz_{ij}^r$, where all $z_{ij}^r$ are independent standard real Gaussian variables. Writing $G_{r,N}:=N^{-1/2}(z_{ij}^r)_{i,j\in\{1,\ldots,N\}}$, we then have
\begin{equation}\label{eq:gaussian-block-pencil}
M_N=\frac12\sum_{r=1}^m\bigl(A_r\otimes G_{r,N}+A_r^{\mathsf T}\otimes G_{r,N}^{\mathsf T}\bigr).
\end{equation}
We use \eqref{eq:gaussian-block-pencil} only under this additional Gaussian assumption or for a Gaussian matrix with the same covariance.

\begin{proposition}\label{prop:strong-convergence}
Under \eqref{eq:coupling-covariance}, there is one probability-one event on which \eqref{eq:affine-spectrum} and \eqref{eq:affine-trace} hold simultaneously for every $B=B^*\in M_d(\C)$, $R\in M_d(\C)$, $t\in\R$, and continuous $f:\R\to\C$. For each fixed $f$, both limits are locally uniform in $(B,R,t)$, and hence hold along every convergent coefficient sequence, including sequences depending on the sample. The covariance map $Q\mapsto(\operatorname{id}\otimes\tau)[\mathsf M(Q\otimes1)\mathsf M]$ of $\mathsf M$ is $\eta$ in \eqref{eq:eta}.
\end{proposition}

\begin{proof}
Appendix~\ref{app:fourth-moment} proves \eqref{eq:affine-spectrum} and \eqref{eq:affine-trace} under \eqref{eq:coupling-covariance}, with the stated simultaneous convergence and local uniformity. In the Gaussian case, \cite[Theorem 3.3 and Remark 3.6]{xiang-zhang} also give joint strong $*$-convergence of $(G_{r,N})_{r=1}^m$ to $(c_r)_{r=1}^m$, including every fixed matrix amplification. This identifies the Gaussian comparison limit in the appendix with \eqref{eq:free-operator}.

The coefficients in \eqref{eq:free-operator-semicircular} are Hermitian, and their covariance map is
\begin{equation*}
\frac18\sum_{r=1}^m\left((A_r+A_r^{\mathsf T})Q(A_r+A_r^{\mathsf T})
-(A_r-A_r^{\mathsf T})Q(A_r-A_r^{\mathsf T})\right)
=\frac14\sum_{r=1}^m(A_rQA_r^{\mathsf T}+A_r^{\mathsf T}QA_r)
=\eta(Q).
\end{equation*}
\end{proof}

\subsection{The matrix Dyson equation}

\begin{proposition}\label{prop:dyson-evaluation}
For $\lambda\in\Lambda_\beta$, the matrix $Q_\lambda$ in \eqref{eq:Q-lambda} satisfies \eqref{eq:matrix-dyson} and
\begin{equation}\label{eq:trace-log-evaluation}
d\tau_d\log(D_\lambda\otimes1-2\beta\mathsf M)
=-\log\det Q_\lambda+2\beta^2\Tr\bigl(Q_\lambda\eta(Q_\lambda)\bigr).
\end{equation}
Consequently, \eqref{eq:main-formula} and \eqref{eq:finite-dimensional-formula} agree.
\end{proposition}

\begin{proof}
For $0\leq u\leq1$, set
\begin{equation*}
B_u:=D_\lambda\otimes1-2\beta u\mathsf M,
\qquad\text{and}\qquad
Q_u:=(\operatorname{id}\otimes\tau)(B_u^{-1}).
\end{equation*}
Since $(\operatorname{id}\otimes\tau)(\mathsf M)=0$, we have $D_\lambda=(\operatorname{id}\otimes\tau)(B_1)\succ0$. The identity
\begin{equation*}
B_u=(1-u)(D_\lambda\otimes1)+uB_1
\end{equation*}
shows that every $B_u$ is strictly positive. By \eqref{eq:free-operator-semicircular}, put
\begin{equation*}
V_{r,1}:=\frac{A_r+A_r^{\mathsf T}}{2\sqrt2},
\qquad
V_{r,2}:=\frac{\mathrm i(A_r-A_r^{\mathsf T})}{2\sqrt2},
\qquad\text{and}\qquad
\widetilde V_{r,k,u}:=2\beta uV_{r,k}.
\end{equation*}
Then $V_{r,k}$ is Hermitian and $\eta(Q)=\sum_{r=1}^m\sum_{k=1}^2V_{r,k}QV_{r,k}$. Moreover,
\begin{equation*}
2\beta u\mathsf M=\sum_{r=1}^m\sum_{k=1}^2\widetilde V_{r,k,u}\otimes s_{r,k}.
\end{equation*}
The semicircular integration-by-parts identity $\tau(s_{r,k}p)=(\tau\otimes\tau)(\partial_{r,k}p)$ gives
\begin{equation*}
(\operatorname{id}\otimes\tau)(2\beta u\mathsf M B_u^{-1})
=\sum_{r=1}^m\sum_{k=1}^2\widetilde V_{r,k,u}Q_u\widetilde V_{r,k,u}Q_u
=4\beta^2u^2\eta(Q_u)Q_u.
\end{equation*}
Indeed, entrywise differentiation gives
\begin{equation*}
\partial_{r,k}(B_u^{-1})_{bc}
=\sum_{p,q=1}^d(B_u^{-1})_{bp}(\widetilde V_{r,k,u})_{pq}\otimes(B_u^{-1})_{qc},
\end{equation*}
so applying $\tau\otimes\tau$ gives $(Q_u\widetilde V_{r,k,u}Q_u)_{bc}$. There are constants $0<\varepsilon<C<\infty$ such that $\operatorname{spec}(B_u)\subset[\varepsilon,C]$ for $0\leq u\leq1$. Polynomial approximation of $t\mapsto t^{-1}$ in $C^1([\varepsilon,C])$ also approximates its divided difference uniformly on $[\varepsilon,C]^2$, and therefore extends the polynomial integration-by-parts identity to $B_u^{-1}$. Applying $\operatorname{id}\otimes\tau$ to $B_uB_u^{-1}=1$ therefore yields
\begin{equation}\label{eq:dyson-path}
Q_u=\left(D_\lambda-4\beta^2u^2\eta(Q_u)\right)^{-1}.
\end{equation}
At $u=0$, $Q_0=D_\lambda^{-1}$; at $u=1$, \eqref{eq:dyson-path} is \eqref{eq:matrix-dyson}.

The map $\eta$ is self-adjoint for the trace inner product. Define
\begin{equation*}
\Psi_u(Q):=-\log\det Q+\Tr(D_\lambda Q)-d-2\beta^2u^2\Tr(Q\eta(Q)).
\end{equation*}
Equation \eqref{eq:dyson-path} is the stationary equation of $Q\mapsto\Psi_u(Q)$. Therefore
\begin{equation}\label{eq:action-derivative}
\frac{\d}{\d u}\Psi_u(Q_u)=-4\beta^2u\Tr(Q_u\eta(Q_u)).
\end{equation}
Tracial differentiation of the logarithm and the Wick identity give
\begin{align*}
\frac{\d}{\d u}d\tau_d\log B_u
&=-2\beta d\tau_d(\mathsf M B_u^{-1})\\
&\stackrel{\text{Wick}}{=}-4\beta^2u\Tr(Q_u\eta(Q_u))
\stackrel{\eqref{eq:action-derivative}}{=}\frac{\d}{\d u}\Psi_u(Q_u).
\end{align*}
At $u=0$, both $d\tau_d\log B_u$ and $\Psi_u(Q_u)$ equal $\log\det D_\lambda$. Integrating and using
\begin{equation*}
\Tr(D_\lambda Q_1)-d=4\beta^2\Tr(Q_1\eta(Q_1))
\end{equation*}
proves \eqref{eq:trace-log-evaluation}.
\end{proof}

\begin{proof}[Proof of Corollary~\ref{cor:finite-dimensional}]
Proposition~\ref{prop:dyson-evaluation} proves \eqref{eq:finite-dimensional-formula}. It remains to identify $\Lambda_\beta$. In \eqref{eq:free-operator-semicircular}, the Hermitian coefficients of the variance-one semicircular variables are
\begin{equation*}
\frac{A_r+A_r^{\mathsf T}}{2\sqrt2}
\quad\text{and}\quad
\frac{\mathrm i(A_r-A_r^{\mathsf T})}{2\sqrt2}.
\end{equation*}
The covariance operator in \cite[Equation (22)]{kunisky} is $\eta$. Apply Kunisky's right-edge formula \cite[Theorem 1.9 and Equation (10)]{kunisky} to $2\beta\mathsf M-D_\lambda\otimes1$. Its covariance map is $4\beta^2\eta$, and hence
\begin{equation}\label{eq:scaled-edge}
\lambda_{\max}(2\beta\mathsf M-D_\lambda\otimes1)
=\max_{\substack{S=S^*,\ S\succeq0\\\Tr S=1}}
\left\{-\Tr(D_\lambda S)+4\beta\Tr\left[(S^\frac{1}{2}\eta(S)S^\frac{1}{2})^\frac{1}{2}\right]\right\}.
\end{equation}
The pencil $D_\lambda\otimes1-2\beta\mathsf M$ is strictly positive exactly when the left-hand side of \eqref{eq:scaled-edge} is negative. This proves \eqref{eq:admissible-finite}.
\end{proof}

\subsection{Rank-one extraction and the determinant reduction}

\begin{lemma}\label{lem:rank-one-extraction}
Under the assumptions of Lemma~\ref{lem:spectral-dual}, let $X_N\in M_{dN}(\R)$ be symmetric and positive semidefinite, and suppose that
\begin{equation*}
\max_{a\in\{1,\ldots,d\}}\frac1N\Tr(P_aX_N)=O(1).
\end{equation*}
Then there are vectors $x_N\in\R^{dN}$ such that
\begin{equation}\label{eq:rank-one-masses}
\|P_ax_N\|^2=\Tr(P_aX_N)
\quad\text{for every }a\in\{1,\ldots,d\},
\end{equation}
and
\begin{equation}\label{eq:rank-one-energy}
x_N^{\mathsf T}M_Nx_N\leq\Tr(M_NX_N)+o(N).
\end{equation}
\end{lemma}

\begin{proof}
Set
\begin{equation*}
R_N:=\diag\left(\sqrt{\frac1N\Tr(P_1X_N)},\ldots,\sqrt{\frac1N\Tr(P_dX_N)}\right).
\end{equation*}
Pass to an arbitrary subsequence and then to a further subsequence such that
\begin{equation*}
\lim_{N\to\infty}\|R_N-\bar R\|=0
\end{equation*}
for some diagonal matrix $\bar R\succeq0$. Define
\begin{align*}
r_N&:=\min\left\{u^{\mathsf T}(R_N\otimes I_N)M_N(R_N\otimes I_N)u:\|P_au\|^2=1\text{ for every }a\in\{1,\ldots,d\}\right\},\\
s_N&:=\min\left\{\Tr((R_N\otimes I_N)M_N(R_N\otimes I_N)Y):Y\succeq0,\ \Tr(P_aY)=1\text{ for every }a\in\{1,\ldots,d\}\right\}.
\end{align*}
The local uniformity established at the start of the proof of Lemma~\ref{lem:spectral-dual} shows that Lemma~\ref{lem:spectral-dual} applies to $-(R_N\otimes I_N)M_N(R_N\otimes I_N)$: its matrices in \eqref{eq:affine-pencils} have coefficients $(B,SR_N,-t)$, which converge to $(B,S\bar R,-t)$ for each fixed $B,S,t$.

Lemma~\ref{lem:spectral-dual}, applied to $-(R_N\otimes I_N)M_N(R_N\otimes I_N)$, gives
\begin{equation}\label{eq:rank-one-limit}
\lim_{N\to\infty}r_N
=\sup_{\theta\in\R^d}\left\{\sum_{a=1}^d\theta_a+d\lambda_{\min}\bigl((\bar R\otimes1)\mathsf M(\bar R\otimes1)-D_\theta\otimes1\bigr)\right\}.
\end{equation}
The point $I_{dN}/N$ satisfies the semidefinite constraints and is positive definite. Strong semidefinite duality, followed by a common scalar shift of $\theta$, gives
\begin{equation}\label{eq:sdp-dual}
s_N=\sup_{\theta\in\R^d}
\left\{\sum_{a=1}^d\theta_a+d\lambda_{\min}\bigl((R_N\otimes I_N)M_N(R_N\otimes I_N)-D_\theta\otimes I_N\bigr)\right\}.
\end{equation}
The right-hand side is invariant under $\theta\mapsto\theta+c\1$. Equation \eqref{eq:affine-spectrum} permits a constant $C_R<\infty$ such that
\begin{equation*}
\sup_N\|(R_N\otimes I_N)M_N(R_N\otimes I_N)\|\leq C_R.
\end{equation*}
On $\sum_a\theta_a=0$, testing a unit vector in a block attaining $\max_a\theta_a$ gives
\begin{equation*}
\sum_{a=1}^d\theta_a
+d\lambda_{\min}\bigl((R_N\otimes I_N)M_N(R_N\otimes I_N)-D_\theta\otimes I_N\bigr)
\leq d\bigl(C_R-\max_a\theta_a\bigr).
\end{equation*}
At $\theta=0$, the objective is at least $-dC_R$. Thus every maximizer satisfies
\begin{equation}\label{eq:theta-compact-bounds}
\max_a\theta_a\leq2C_R,
\qquad\text{and}\qquad
\min_a\theta_a\geq-2(d-1)C_R.
\end{equation}
Equation \eqref{eq:affine-spectrum}, the $\|\cdot\|_\infty$-Lipschitz dependence on $\theta$, and a finite net on the compact set in \eqref{eq:theta-compact-bounds} show that the limit of \eqref{eq:sdp-dual} equals the right-hand side of \eqref{eq:rank-one-limit}. Therefore
\begin{equation}\label{eq:rank-sdp-gap}
\lim_{N\to\infty}(r_N-s_N)=0,
\end{equation}
along the chosen subsequence. Since the initial subsequence was arbitrary, \eqref{eq:rank-sdp-gap} holds for the full sequence; otherwise a subsequence bounded away from zero would have a further subsequence on which \eqref{eq:rank-sdp-gap} holds.

There exists $Y_N\succeq0$ with $\Tr(P_aY_N)=1$ for every $a\in\{1,\ldots,d\}$ such that
\begin{equation*}
\frac{X_N}{N}=(R_N\otimes I_N)Y_N(R_N\otimes I_N).
\end{equation*}
Indeed, if $\Tr(P_aX_N)=0$, positivity forces the corresponding block row and column of $X_N$ to vanish, and the corresponding diagonal block of $Y_N$ may be set equal to $I_N/N$. Hence
\begin{equation*}
s_N\leq\frac1N\Tr(M_NX_N).
\end{equation*}
Choose a minimizer in the definition of $r_N$ and let $x_N$ be $\sqrt N(R_N\otimes I_N)$ times this minimizer. Then \eqref{eq:rank-one-masses} holds, and
\begin{equation*}
x_N^{\mathsf T}M_Nx_N
=Nr_N
\stackrel{\eqref{eq:rank-sdp-gap}}{\leq}Ns_N+o(N)
\leq\Tr(M_NX_N)+o(N).
\end{equation*}
This proves \eqref{eq:rank-one-energy}.
\end{proof}

For $\beta>0$ and $\lambda\in\R^d$, define
\begin{equation}\label{eq:finite-precision-objective}
B_{N,\beta}(\lambda):=D_\lambda\otimes I_N-2\beta M_N
\qquad\text{and}\qquad
\phi_{N,\beta}(\lambda):=
\frac12\sum_{a=1}^d(\lambda_a-1)
-\frac1{2N}\log\det B_{N,\beta}(\lambda),
\end{equation}
whenever $B_{N,\beta}(\lambda)\succ0$.

\begin{proposition}\label{prop:determinant-principle}
Let $M_N$ be as in \eqref{eq:block-pencil} and satisfy the assumptions of Lemma~\ref{lem:spectral-dual}. For every fixed $\beta>0$, we have
\begin{equation}\label{eq:finite-determinant-principle}
\lim_{N\to\infty}\left[
F_N(\beta)-\inf_{\substack{\lambda\in\R^d\\B_{N,\beta}(\lambda)\succ0}}\phi_{N,\beta}(\lambda)
\right]=0.
\end{equation}
\end{proposition}

\begin{proof}
\smallskip
\noindent\emph{Step 1: Existence and bounds for $\lambda_{N,\beta}$.}
Set
\begin{equation*}
\Phi_{N,\beta}:=\inf_{\substack{\lambda\in\R^d\\B_{N,\beta}(\lambda)\succ0}}\phi_{N,\beta}(\lambda).
\end{equation*}
Its domain contains every $c\1$ with $c>2\beta\|M_N\|$. Positivity of the $a$th diagonal principal block gives
\begin{equation}\label{eq:lambda-lower-finite}
\lambda_a>-2\beta\|M_N\|.
\end{equation}
Moreover,
\begin{multline*}
\frac12\sum_{a=1}^d(\lambda_a-1)
-\frac1{2N}\log\det(D_\lambda\otimes I_N-2\beta M_N)\\
\geq\frac12\left(\max_a\lambda_a-2\beta(d-1)\|M_N\|-d\right)
-\frac d2\log\left(\max_a\lambda_a+2\beta\|M_N\|\right),
\end{multline*}
and the right-hand side tends to $+\infty$ as $\max_a\lambda_a\to\infty$. If $(\lambda^{(k)})$ lies in the domain and converges to a finite boundary point, then, with $B^{(k)}:=D_{\lambda^{(k)}}\otimes I_N-2\beta M_N$,
\begin{equation*}
\sup_k\|B^{(k)}\|<\infty,
\qquad
\lim_{k\to\infty}\lambda_{\min}(B^{(k)})=0,
\qquad\text{and}\qquad
\lim_{k\to\infty}\bigl[-\log\det B^{(k)}\bigr]=+\infty.
\end{equation*}
Finally, if $B_{N,\beta}(\lambda)\succ0$ and $v\ne0$, then
\begin{equation}\label{eq:finite-hessian}
\begin{aligned}
&\sum_{a,b=1}^dv_av_b\frac{\partial^2}{\partial\lambda_a\partial\lambda_b}
\left[-\frac1{2N}\log\det B_{N,\beta}(\lambda)\right]\\
&\qquad=\frac1{2N}\Tr\left[
\left(B_{N,\beta}(\lambda)^{-1/2}(D_v\otimes I_N)B_{N,\beta}(\lambda)^{-1/2}\right)^2
\right]>0.
\end{aligned}
\end{equation}
Equations \eqref{eq:lambda-lower-finite}--\eqref{eq:finite-hessian} prove existence and uniqueness of a minimizer $\lambda_{N,\beta}$. Differentiation at the minimizer gives
\begin{equation}\label{eq:inverse-block-masses}
\Tr(P_aB_{N,\beta}(\lambda_{N,\beta})^{-1})=N
\quad\text{for every }a\in\{1,\ldots,d\}.
\end{equation}
Taking the trace of $P_aB_{N,\beta}(\lambda_{N,\beta})B_{N,\beta}(\lambda_{N,\beta})^{-1}=P_a$ yields
\begin{equation}\label{eq:lambda-stationarity}
(\lambda_{N,\beta})_a=1+\frac{2\beta}{N}\Tr(P_aM_NB_{N,\beta}(\lambda_{N,\beta})^{-1}).
\end{equation}
By Hilbert--Schmidt Cauchy--Schwarz and \eqref{eq:inverse-block-masses},
\begin{equation*}
|\Tr(P_aM_NB_{N,\beta}(\lambda_{N,\beta})^{-1})|
\leq\bigl(\Tr(P_aB_{N,\beta}(\lambda_{N,\beta})^{-1})\bigr)^\frac{1}{2}
\bigl(\Tr(B_{N,\beta}(\lambda_{N,\beta})^{-1}M_N^2)\bigr)^\frac{1}{2}
\leq\sqrt d\,N\|M_N\|.
\end{equation*}
Consequently, we have
\begin{equation*}
\sup_N\|\lambda_{N,\beta}\|<\infty,
\qquad\text{and}\qquad
\sup_N\|B_{N,\beta}(\lambda_{N,\beta})\|<\infty.
\end{equation*}

\smallskip
\noindent\emph{Step 2: Upper bound in \eqref{eq:finite-determinant-principle}.}
We first prove the upper bound. The Gibbs variational principle on $\mathcal S_N$ gives
\begin{equation}\label{eq:gibbs-sphere}
NF_N(\beta)=\sup_\nu
\left\{\beta\E_\nu x^{\mathsf T}M_Nx-D(\nu\|\mu_N^{\otimes d})\right\}.
\end{equation}
Given $x\sim\nu$, independently multiply each block $x^a$ by $N^{-1/2}$ times the Euclidean norm of a standard Gaussian vector in $\R^N$, and denote the resulting vector by $X$. The polar decomposition of standard Gaussian measure gives
\begin{equation}\label{eq:radial-entropy}
D(\operatorname{Law}(X)\|\gamma_{dN})=D(\nu\|\mu_N^{\otimes d}),
\qquad\text{and}\qquad
\E\|P_aX\|^2=N,
\end{equation}
where $\gamma_{dN}$ is standard Gaussian measure. Moreover,
\begin{equation}\label{eq:radial-energy-error}
\begin{aligned}
\E|X^{\mathsf T}M_NX-x^{\mathsf T}M_Nx|
&\leq\|M_N\|\E\bigl[\|X-x\|(\|X\|+\|x\|)\bigr]\\
&\stackrel{\text{Cauchy--Schwarz}}{\leq}\|M_N\|
\bigl(\E\|X-x\|^2\bigr)^\frac{1}{2}
\bigl(\E(\|X\|+\|x\|)^2\bigr)^\frac{1}{2}
=O(\sqrt N).
\end{aligned}
\end{equation}
Here $\E\|X-x\|^2=O(1)$ and $\E(\|X\|+\|x\|)^2=O(N)$.

For any $\lambda$ such that $B_{N,\beta}(\lambda)\succ0$, the Gaussian entropy variational principle and
\begin{equation*}
\log\int_{\R^{dN}}
\exp\left(\frac12\|z\|^2-\frac12z^{\mathsf T}B_{N,\beta}(\lambda)z\right)\,\d\gamma_{dN}(z)
=-\frac12\log\det B_{N,\beta}(\lambda)
\end{equation*}
gives
\begin{equation}\label{eq:gaussian-upper}
-\frac12\E X^{\mathsf T}B_{N,\beta}(\lambda)X-D(\operatorname{Law}(X)\|\gamma_{dN})
\leq-\frac{dN}2-\frac12\log\det B_{N,\beta}(\lambda).
\end{equation}
Combining \eqref{eq:gibbs-sphere}, \eqref{eq:radial-entropy}, \eqref{eq:radial-energy-error}, and \eqref{eq:gaussian-upper}, and then minimizing over $\lambda$, yields
\begin{equation}\label{eq:determinant-upper}
F_N(\beta)\leq\Phi_{N,\beta}+O(N^{-1/2}).
\end{equation}

\smallskip
\noindent\emph{Step 3: Construction of the Gaussian trial measure.}
For the lower bound, fix $\ell>1$ and use functional calculus with $B_{N,\beta}(\lambda_{N,\beta})$ to define
\begin{equation}\label{eq:spectral-cap}
\Sigma_{N,\ell}:=\min\{B_{N,\beta}(\lambda_{N,\beta})^{-1},\ell I_{dN}\},
\qquad\text{and}\qquad
K_{N,\ell}:=B_{N,\beta}(\lambda_{N,\beta})^{-1}-\Sigma_{N,\ell}\succeq0,
\end{equation}
where the minimum is defined by functional calculus. Equation \eqref{eq:inverse-block-masses} gives
\begin{equation}\label{eq:deficit-masses}
\Tr(P_aK_{N,\ell})=N-\Tr(P_a\Sigma_{N,\ell}).
\end{equation}
On an arbitrary subsequence, pass to a further subsequence along which $(\lambda_{N,\beta})$ converges. For every fixed $C=C^*\in M_d(\C)$, $R\in M_d(\C)$, and $t\in\R$, we have
\begin{equation*}
C\otimes I_N+t(R\otimes I_N)B_{N,\beta}(\lambda_{N,\beta})(R^*\otimes I_N)
=L_N(C+tRD_{\lambda_{N,\beta}}R^*,R,-2\beta t).
\end{equation*}
The local uniformity in Lemma~\ref{lem:spectral-dual} therefore gives \eqref{eq:affine-spectrum}--\eqref{eq:affine-trace} for $B_{N,\beta}(\lambda_{N,\beta})$. Apply Lemma~\ref{lem:rank-one-extraction} with $B_{N,\beta}(\lambda_{N,\beta})$ in place of $M_N$ and $K_{N,\ell}$ in place of $X_N$. It gives a vector $m_{N,\ell}$ satisfying
\begin{equation}\label{eq:mean-masses}
\|P_am_{N,\ell}\|^2=\Tr(P_aK_{N,\ell})
\quad\text{and}\quad
m_{N,\ell}^{\mathsf T}B_{N,\beta}(\lambda_{N,\beta})m_{N,\ell}\leq\Tr(B_{N,\beta}(\lambda_{N,\beta})K_{N,\ell})+o(N).
\end{equation}
Since $-\log(\ell b)\leq(\ell b)^{-1}$ for $0<b<\ell^{-1}$, functional calculus, \eqref{eq:spectral-cap}, and \eqref{eq:inverse-block-masses} give
\begin{equation}\label{eq:cap-log-loss}
0\leq-\log\det(B_{N,\beta}(\lambda_{N,\beta})\Sigma_{N,\ell})
\leq\frac1\ell\Tr(B_{N,\beta}(\lambda_{N,\beta})^{-1})
=\frac{dN}{\ell}.
\end{equation}
Let $Z_{N,\ell}$ be Gaussian with mean $m_{N,\ell}$ and covariance $\Sigma_{N,\ell}$. By \eqref{eq:deficit-masses} and \eqref{eq:mean-masses},
\begin{equation}\label{eq:gaussian-spherical-moments}
\E\|P_aZ_{N,\ell}\|^2=N
\quad\text{for every }a\in\{1,\ldots,d\}.
\end{equation}
Define $Y_{N,\ell}\in\mathcal S_N$ by
\begin{equation*}
P_aY_{N,\ell}:=\frac{\sqrt N\,P_aZ_{N,\ell}}{\|P_aZ_{N,\ell}\|}
\quad\text{for every }a\in\{1,\ldots,d\}.
\end{equation*}
For fixed $\ell$, \eqref{eq:spectral-cap} and $\sup_N\|B_{N,\beta}(\lambda_{N,\beta})\|<\infty$ bound $\Sigma_{N,\ell}$ above and below by positive scalar multiples of $I_{dN}$, uniformly in $N$. The noncentral Gaussian variance formula and \eqref{eq:gaussian-spherical-moments} give
\begin{equation}\label{eq:block-norm-variance}
\operatorname{Var}(\|P_aZ_{N,\ell}\|^2)
=2\Tr(P_a\Sigma_{N,\ell}P_a\Sigma_{N,\ell})
+4m_{N,\ell}^{\mathsf T}P_a\Sigma_{N,\ell}P_am_{N,\ell}
=O_{\ell}(N).
\end{equation}
Since
\begin{equation*}
(\|P_aZ_{N,\ell}\|-\sqrt N)^2\leq N^{-1}(\|P_aZ_{N,\ell}\|^2-N)^2,
\end{equation*}
equation \eqref{eq:block-norm-variance} implies $\E\|Y_{N,\ell}-Z_{N,\ell}\|^2=O_{\ell}(1)$. Hence
\begin{equation}\label{eq:projected-energy}
\begin{aligned}
\left|\E Y_{N,\ell}^{\mathsf T}M_NY_{N,\ell}-\E Z_{N,\ell}^{\mathsf T}M_NZ_{N,\ell}\right|
&\stackrel{\text{Cauchy--Schwarz}}{\leq}\|M_N\|
\bigl(\E\|Y_{N,\ell}-Z_{N,\ell}\|^2\bigr)^\frac{1}{2}\\
&\qquad\times\bigl(\E(\|Y_{N,\ell}\|+\|Z_{N,\ell}\|)^2\bigr)^\frac{1}{2}
=O_{\ell}(\sqrt N),\\
\E Z_{N,\ell}^{\mathsf T}M_NZ_{N,\ell}
&=m_{N,\ell}^{\mathsf T}M_Nm_{N,\ell}+\Tr(M_N\Sigma_{N,\ell}).
\end{aligned}
\end{equation}

\smallskip
\noindent\emph{Step 4: Entropy and the limit $\ell\to\infty$.}
Block normalization sends $\gamma_{dN}$ to $\mu_N^{\otimes d}$. The data-processing inequality and the Gaussian relative-entropy formula therefore give
\begin{align}
D(\operatorname{Law}(Y_{N,\ell})\|\mu_N^{\otimes d})
&\leq D(\mathcal N(m_{N,\ell},\Sigma_{N,\ell})\|\mathcal N(0,I_{dN}))\notag\\
&=\frac12\left(\|m_{N,\ell}\|^2+\Tr\Sigma_{N,\ell}-dN-\log\det\Sigma_{N,\ell}\right)\notag\\
&\stackrel{\eqref{eq:gaussian-spherical-moments}}{=}-\frac12\log\det\Sigma_{N,\ell}.
\label{eq:projected-entropy}
\end{align}
Use the law of $Y_{N,\ell}$ in \eqref{eq:gibbs-sphere}. Equations \eqref{eq:projected-energy} and \eqref{eq:projected-entropy} give
\begin{equation*}
F_N(\beta)\geq
\frac\beta N\left(m_{N,\ell}^{\mathsf T}M_Nm_{N,\ell}+\Tr(M_N\Sigma_{N,\ell})\right)
+\frac1{2N}\log\det\Sigma_{N,\ell}+o(1).
\end{equation*}
Subtracting this lower bound from $\Phi_{N,\beta}$ and using $2\beta M_N=D_{\lambda_{N,\beta}}\otimes I_N-B_{N,\beta}(\lambda_{N,\beta})$, we obtain
\begin{align*}
&\Phi_{N,\beta}-\frac\beta N\left(m_{N,\ell}^{\mathsf T}M_Nm_{N,\ell}+\Tr(M_N\Sigma_{N,\ell})\right)-\frac1{2N}\log\det\Sigma_{N,\ell}\\
&\qquad=\frac{m_{N,\ell}^{\mathsf T}B_{N,\beta}(\lambda_{N,\beta})m_{N,\ell}+\Tr(B_{N,\beta}(\lambda_{N,\beta})\Sigma_{N,\ell})-dN-\log\det(B_{N,\beta}(\lambda_{N,\beta})\Sigma_{N,\ell})}{2N}\\
&\qquad\stackrel{\eqref{eq:mean-masses},\eqref{eq:spectral-cap}}{\leq}
-\frac1{2N}\log\det(B_{N,\beta}(\lambda_{N,\beta})\Sigma_{N,\ell})+o(1)
\stackrel{\eqref{eq:cap-log-loss}}{\leq}\frac d{2\ell}+o(1).
\end{align*}
Thus, along the chosen further subsequence,
\begin{equation}\label{eq:determinant-lower-cutoff}
\liminf_{N\to\infty}\bigl(F_N(\beta)-\Phi_{N,\beta}\bigr)\geq-\frac d{2\ell}.
\end{equation}
Since the initial subsequence was arbitrary, \eqref{eq:determinant-lower-cutoff} holds for the full sequence. Letting $\ell\to\infty$ and combining with \eqref{eq:determinant-upper} proves \eqref{eq:finite-determinant-principle}.
\end{proof}

\subsection{Proofs of the limiting formulas}

\begin{proof}[Proof of Theorem~\ref{thm:main}]
Work on the probability-one event in Proposition~\ref{prop:strong-convergence}. Proposition~\ref{prop:strong-convergence} verifies \eqref{eq:affine-spectrum}--\eqref{eq:affine-trace}, which are the hypotheses of Lemma~\ref{lem:spectral-dual} and hence of Proposition~\ref{prop:determinant-principle}. Therefore Proposition~\ref{prop:determinant-principle} reduces $F_N(\beta)$ to the infimum of $\phi_{N,\beta}$ over the domain where $B_{N,\beta}(\lambda)\succ0$. Recall $\phi_\beta$ from \eqref{eq:limiting-objective}. For fixed $\lambda\in\Lambda_\beta$, \eqref{eq:affine-spectrum} places the finite spectra in a common compact subinterval of $(0,\infty)$ for all sufficiently large $N$. Applying \eqref{eq:affine-trace} to a continuous extension of the logarithm from this interval gives
\begin{equation}\label{eq:fixed-lambda-log-limit}
\lim_{N\to\infty}\phi_{N,\beta}(\lambda)=\phi_\beta(\lambda).
\end{equation}
By \eqref{eq:fixed-lambda-log-limit},
\begin{equation}\label{eq:inf-limsup}
\limsup_{N\to\infty}\inf\phi_{N,\beta}\leq\inf_{\lambda\in\Lambda_\beta}\phi_\beta(\lambda).
\end{equation}

Let $\lambda_{N,\beta}$ be the unique minimizer from the proof of Proposition~\ref{prop:determinant-principle}. Equation \eqref{eq:lambda-stationarity} shows that $\lambda_{N,\beta}$ is bounded. Along any subsequence, pass to a further subsequence and choose $\lambda$ such that
\begin{equation*}
\lim_{N\to\infty}\|\lambda_{N,\beta}-\lambda\|=0.
\end{equation*}
Equation \eqref{eq:affine-spectrum} and local uniformity then give
\begin{equation*}
D_\lambda\otimes1-2\beta\mathsf M\succeq0.
\end{equation*}
For every $\delta>0$, one has $\lambda_{N,\beta}\leq\lambda+\delta\1$ componentwise for all sufficiently large $N$. Monotonicity of the logarithm and \eqref{eq:affine-trace} at the strictly positive pencil with parameter $\lambda+\delta\1$ give
\begin{align}
\liminf_{N\to\infty}\phi_{N,\beta}(\lambda_{N,\beta})
&\stackrel{\text{mono.}}{\geq}\frac12\sum_{a=1}^d(\lambda_a-1)
-\frac d2\tau_d\log(D_\lambda\otimes1-2\beta\mathsf M+\delta\1)\stackrel{\eqref{eq:limiting-objective}}{=}\phi_\beta(\lambda+\delta\1)-\frac{d\delta}2\notag\\
&\geq\inf_{\mu\in\Lambda_\beta}\phi_\beta(\mu)-\frac{d\delta}2.
\label{eq:inf-liminf}
\end{align}
Let $\delta\downarrow0$. Since the initial subsequence was arbitrary, \eqref{eq:inf-limsup} and \eqref{eq:inf-liminf} prove \eqref{eq:main-formula} almost surely.

Finally, assume that $g_{11}$ is Gaussian and use \eqref{eq:gaussian-block-pencil}. Then
\begin{equation}\label{eq:free-energy-norm-bound}
|F_N(\beta)|\leq\beta d\|M_N\|
\leq\beta d\sum_{r=1}^m\|A_r\|\,\|G_{r,N}\|.
\end{equation}
Let $\mathcal E_N$ be a $1/4$-net of the unit sphere with $|\mathcal E_N|\leq9^N$. Since
\begin{equation*}
\|G_{r,N}\|\leq2\max_{u,v\in\mathcal E_N}|u^{\mathsf T}G_{r,N}v|,
\end{equation*}
a union bound for centered Gaussian variables of variance $N^{-1}$ gives, for all sufficiently large $t$,
\begin{equation*}
\P(\|G_{r,N}\|>t)\leq2\cdot81^N e^{-Nt^2/8}.
\end{equation*}
Consequently, for every $q<\infty$, we have
\begin{equation}\label{eq:ginibre-norm-moments}
\sup_N\E\|G_{r,N}\|^q<\infty.
\end{equation}
For fixed $p<q$, \eqref{eq:free-energy-norm-bound} therefore gives
\begin{equation*}
\sup_N\E|F_N(\beta)|^q<\infty.
\end{equation*}
Hence $|F_N(\beta)-F(\beta)|^p$ is uniformly integrable, and almost-sure convergence implies convergence in $L^p$.
\end{proof}

\begin{lemma}\label{lem:fidelity-concavity}
The function
\begin{equation*}
S\longmapsto\Tr\left[(S^\frac{1}{2}\eta(S)S^\frac{1}{2})^\frac{1}{2}\right]
\end{equation*}
is concave and continuous on the positive semidefinite cone in $M_d(\C)$.
\end{lemma}

\begin{proof}
For positive semidefinite matrices $B$ and $C$, we have
\begin{equation}\label{eq:fidelity-sdp}
2\Tr[(B^\frac{1}{2}CB^\frac{1}{2})^\frac{1}{2}]
=\max\left\{2\operatorname{Re}\Tr X:
\begin{pmatrix}B&X\\X^*&C\end{pmatrix}\succeq0\right\}.
\end{equation}
Indeed, block positivity is equivalent to $X=B^\frac{1}{2}KC^\frac{1}{2}$ for a contraction $K$, and the polar decomposition evaluates the maximum. Apply \eqref{eq:fidelity-sdp} with $B=S$ and $C=\eta(S)$. Since $\eta$ is linear, convex combinations of block matrices satisfying the positivity constraint in \eqref{eq:fidelity-sdp} also satisfy that constraint, which proves concavity. Continuity follows from continuity of the square root and trace norm.
\end{proof}

\begin{proof}[Proof of Theorem~\ref{thm:ground-state}]
\smallskip
\noindent\emph{Step 1: The limit $\beta\to\infty$.}
With $\phi_\beta$ defined by \eqref{eq:limiting-objective}, set
\begin{equation}\label{eq:positive-dual}
\mathcal G:=\inf\left\{\sum_{a=1}^d\mu_a:D_\mu\otimes1-\mathsf M\succeq0\right\}.
\end{equation}
If $D_\mu\otimes1-\mathsf M\succeq0$, then $\lambda:=2\beta(\mu+\varepsilon\1)$ belongs to $\Lambda_\beta$ for every $\varepsilon>0$. Substitution in \eqref{eq:main-formula} gives
\begin{align*}
\frac{F(\beta)}\beta
\leq\sum_{a=1}^d\mu_a+d\varepsilon-\frac d{2\beta}
-\frac d{2\beta}\log(2\beta)
-\frac d{2\beta}\tau_d\log(D_\mu\otimes1-\mathsf M+\varepsilon\1),
\end{align*}
and therefore
\begin{equation*}
\limsup_{\beta\to\infty}\frac{F(\beta)}\beta\leq\sum_{a=1}^d\mu_a+d\varepsilon.
\end{equation*}
Taking first the infimum over $\mu$ and then $\varepsilon\downarrow0$ proves the upper bound by $\mathcal G$.

For the reverse bound, choose $\widetilde\lambda_\beta\in\Lambda_\beta$ such that $\phi_\beta(\widetilde\lambda_\beta)\leq F(\beta)+1$, and set $c:=\|\mathsf M\|+1$. Since $2\beta c\1\in\Lambda_\beta$, the choice of $\widetilde\lambda_\beta$ gives
\begin{equation*}
\phi_\beta(\widetilde\lambda_\beta)\leq F(\beta)+1\leq\phi_\beta(2\beta c\1)+1.
\end{equation*}
Moreover,
\begin{equation*}
\phi_\beta(2\beta c\1)
=\beta dc-\frac d2-\frac d2\log(2\beta)
-\frac d2\tau_d\log(c\1-\mathsf M)
\leq C(1+\beta),
\end{equation*}
for a constant $C$ independent of $\beta\geq1$. Positivity of the $a$th diagonal compression gives
\begin{equation}\label{eq:zero-temperature-component-lower}
(\widetilde\lambda_\beta)_a>-2\beta\|\mathsf M\|.
\end{equation}
We have $\|D_{\widetilde\lambda_\beta}\otimes1-2\beta\mathsf M\|\leq\max_a(\widetilde\lambda_\beta)_a+2\beta\|\mathsf M\|$ and
\begin{equation*}
\phi_\beta(\widetilde\lambda_\beta)
\geq\frac12\left(\max_a(\widetilde\lambda_\beta)_a-2\beta(d-1)\|\mathsf M\|-d\right)
-\frac d2\log\left(\max_a(\widetilde\lambda_\beta)_a+2\beta\|\mathsf M\|\right),
\end{equation*}
so $\max_a(\widetilde\lambda_\beta)_a=O(\beta)$. Together with \eqref{eq:zero-temperature-component-lower}, this makes $\mu_\beta:=\widetilde\lambda_\beta/(2\beta)$ bounded. For any sequence $(\beta_n)$ satisfying $\lim_{n\to\infty}\beta_n=+\infty$, pass to a subsequence such that
\begin{equation*}
\lim_{n\to\infty}\mu_{\beta_n}=\mu.
\end{equation*}
Since $D_{\mu_{\beta_n}}\otimes1-\mathsf M\succ0$ for every $n$, passing to the limit gives
\begin{equation*}
D_\mu\otimes1-\mathsf M\succeq0.
\end{equation*}
Moreover, we have
\begin{align*}
\frac{\phi_\beta(\widetilde\lambda_\beta)}\beta
&=\sum_{a=1}^d(\mu_\beta)_a-\frac d{2\beta}
-\frac d{2\beta}\tau_d\log(D_{\widetilde\lambda_\beta}\otimes1-2\beta\mathsf M)\\
&\geq\sum_{a=1}^d(\mu_\beta)_a
-\frac d{2\beta}
-\frac d{2\beta}\log\|D_{\widetilde\lambda_\beta}\otimes1-2\beta\mathsf M\|=\sum_{a=1}^d(\mu_\beta)_a-o(1),
\end{align*}
because the norm inside the logarithm is $O(\beta)$. Since $F(\beta)\geq\phi_\beta(\widetilde\lambda_\beta)-1$, the initial sequence $(\beta_n)$ was arbitrary, and every resulting $\mu$ satisfies the constraint $D_\mu\otimes1-\mathsf M\succeq0$ in \eqref{eq:positive-dual}. Hence
\begin{equation}\label{eq:temperature-to-positive-dual}
\lim_{\beta\to\infty}\frac{F(\beta)}\beta=\mathcal G.
\end{equation}

\smallskip
\noindent\emph{Step 2: Identification with the finite-$N$ ground state.}
For $\lambda\in\R^d$, put $t_\lambda:=\lambda_{\max}(\mathsf M-D_\lambda\otimes1)$. Then $D_{\lambda+t_\lambda\1}\otimes1-\mathsf M\succeq0$ and
\begin{equation*}
\sum_{a=1}^d(\lambda_a+t_\lambda)=\sum_{a=1}^d\lambda_a+dt_\lambda.
\end{equation*}
Conversely, if $D_\mu\otimes1-\mathsf M\succeq0$, then $t_\mu\leq0$. These two implications give
\begin{equation}\label{eq:positive-spectral-dual}
\mathcal G=\inf_{\lambda\in\R^d}
\left\{\sum_{a=1}^d\lambda_a+d\lambda_{\max}(\mathsf M-D_\lambda\otimes1)\right\}.
\end{equation}
Lemma~\ref{lem:spectral-dual} and Proposition~\ref{prop:strong-convergence} therefore give
\begin{equation}\label{eq:ground-state-to-positive-dual}
\lim_{N\to\infty}\frac1N\max_{\sigma\in\mathcal S_N}H_N(\sigma)
\stackrel{\eqref{eq:spectral-dual}}{=}\mathcal G.
\end{equation}
This proves the almost-sure assertion under \eqref{eq:coupling-covariance}. If $g_{11}$ is Gaussian, then, for every $q<\infty$, \eqref{eq:ginibre-norm-moments} and the bound
\begin{equation*}
\left|\frac1N\max_{\sigma\in\mathcal S_N}H_N(\sigma)\right|
\leq d\|M_N\|
\leq d\sum_{r=1}^m\|A_r\|\,\|G_{r,N}\|
\end{equation*}
give a uniform $L^q$ bound. Hence \eqref{eq:ground-state-to-positive-dual} also holds in $L^p$ for every $1\leq p<\infty$.

\smallskip
\noindent\emph{Step 3: Evaluation of $\mathcal G$.}
Kunisky's formula \cite[Theorem 1.9 and Equation (10)]{kunisky}, with covariance map \eqref{eq:eta}, gives
\begin{equation}\label{eq:kunisky-ground-state-edge}
\lambda_{\max}(\mathsf M-D_\lambda\otimes1)
=\max_{\substack{S=S^*,\ S\succeq0\\\Tr S=1}}
\left\{-\Tr(D_\lambda S)+2\Tr[(S^\frac{1}{2}\eta(S)S^\frac{1}{2})^\frac{1}{2}]\right\}.
\end{equation}
The density matrices form a compact convex set, and Lemma~\ref{lem:fidelity-concavity} implies that the maximand in \eqref{eq:kunisky-ground-state-edge} is concave in $S$ and affine in $\lambda$. Sion's theorem \cite[Theorem 4.2]{sion} applies with $S$ as the compact variable. Moreover,
\begin{equation*}
\inf_{\lambda\in\R^d}\sum_{a=1}^d\lambda_a(1-dS_{aa})
=\begin{cases}
0,&S_{aa}=d^{-1}\text{ for every }a\in\{1,\ldots,d\},\\
-\infty,&\text{otherwise}.
\end{cases}
\end{equation*}
Therefore \eqref{eq:positive-spectral-dual} becomes
\begin{equation}\label{eq:balanced-density-ground-state}
\mathcal G=2d\max\left\{
\Tr[(S^\frac{1}{2}\eta(S)S^\frac{1}{2})^\frac{1}{2}]:
S=S^*,\ S\succeq0,\ S_{aa}=d^{-1}\text{ for every }a\in\{1,\ldots,d\}
\right\}.
\end{equation}
For complex Hermitian $S$, the real coefficients in \eqref{eq:eta} give
\begin{equation*}
\eta(\bar S)=\overline{\eta(S)},
\qquad\text{and}\qquad
\Tr[(\bar S^\frac{1}{2}\eta(\bar S)\bar S^\frac{1}{2})^\frac{1}{2}]
=\Tr[(S^\frac{1}{2}\eta(S)S^\frac{1}{2})^\frac{1}{2}].
\end{equation*}
The matrix $(S+\bar S)/2$ is real symmetric, positive semidefinite, and has the same diagonal as $S$. Lemma~\ref{lem:fidelity-concavity} gives
\begin{equation*}
\Tr\left[\left(\left(\frac{S+\bar S}{2}\right)^\frac{1}{2}
\eta\left(\frac{S+\bar S}{2}\right)
\left(\frac{S+\bar S}{2}\right)^\frac{1}{2}\right)^\frac{1}{2}\right]
\geq
\Tr\left[(S^\frac{1}{2}\eta(S)S^\frac{1}{2})^\frac{1}{2}\right].
\end{equation*}
Thus the maximum in \eqref{eq:balanced-density-ground-state} may be restricted to real symmetric matrices. Setting $Q:=dS$ and using the linearity of $\eta$, together with \eqref{eq:temperature-to-positive-dual} and \eqref{eq:ground-state-to-positive-dual}, proves \eqref{eq:ground-state-limit}.
\end{proof}

\section{Fluctuations}\label{sec:fluctuations}

Recall $M_N$ from \eqref{eq:block-pencil}, the projections $P_a=E_{aa}\otimes I_N$, and $B_{N,\beta}$ and $\phi_{N,\beta}$ from \eqref{eq:finite-precision-objective}. Recall also $g_\beta$, $\Lambda_\beta$, $\overline{\Lambda}_\beta$, and $\partial\Lambda_\beta$ from \eqref{eq:edge-function}--\eqref{eq:admissible-closure-boundary}, and $\phi_\beta$, $\widehat\phi_\beta$, and $\theta_\beta$ from \eqref{eq:limiting-objective}, \eqref{eq:two-sided-objective}, and Definition~\ref{def:boundary-multiplier}. The Gaussian assumption is imposed only for Corollary~\ref{cor:gaussian-interior-fluctuations} and Subsection~\ref{sec:bipartite-fluctuations}; the other results use \eqref{eq:coupling-covariance} and their stated additional hypotheses. For $\lambda\in\R^d$, define
\begin{equation*}
B_\beta(\lambda):=D_\lambda\otimes1-2\beta\mathsf M
\qquad\text{and}\qquad
g_{N,\beta}(\lambda):=\lambda_{\max}(2\beta M_N-D_\lambda\otimes I_N).
\end{equation*}
We write $\operatorname{TW}_1$ for the GOE Tracy--Widom distribution. For deterministic $a_N>0$, we write $X_N=O_{\mathbb P}(a_N)$ if $(X_N/a_N)_{N\in\mathbb N}$ is tight, and $X_N=o_{\mathbb P}(a_N)$ if $\lim_{N\to\infty}X_N/a_N=0$ in probability.

\subsection{Fluctuations when \texorpdfstring{$\lambda_\beta\in\Lambda_\beta$}{lambda beta is in Lambda beta}}

We first transfer fluctuations of a logarithmic linear statistic to $F_N(\beta)$ under \eqref{eq:coupling-covariance} and an explicit variance bound. We then verify this bound and a Gaussian limit when $g_{11}$ is Gaussian. We begin with uniqueness of a minimizer with finite value in $\overline\Lambda_\beta$.

\begin{lemma}\label{lem:interior-minimizer-uniqueness}
For every $\beta>0$, the function $\phi_\beta$ is strictly convex on $\Lambda_\beta$. There is at most one $\lambda\in\overline\Lambda_\beta$ with finite value $\widehat\phi_\beta(\lambda)=F(\beta)$.
\end{lemma}

\begin{proof}
Fix $\lambda\in\Lambda_\beta$. For every $v\in\R^d\setminus\{0\}$, we have
\begin{equation}\label{eq:limiting-hessian}
\begin{aligned}
v^{\mathsf T}\nabla^2\phi_\beta(\lambda)v
&=\frac d2\tau_d
\left[B_\beta(\lambda)^{-1}(D_v\otimes1)B_\beta(\lambda)^{-1}(D_v\otimes1)\right]\\
&=\frac d2\tau_d
\left[\left(B_\beta(\lambda)^{-1/2}(D_v\otimes1)B_\beta(\lambda)^{-1/2}\right)^2\right]
\geq\frac{\|v\|^2}{2\|B_\beta(\lambda)\|^2}>0.
\end{aligned}
\end{equation}
For the inequality, write $\|X\|_{2,\tau_d}:=\tau_d(X^*X)^\frac{1}{2}$ and apply $\|AXB\|_{2,\tau_d}\leq\|A\|\,\|B\|\,\|X\|_{2,\tau_d}$ to
\begin{equation*}
D_v\otimes1
=B_\beta(\lambda)^\frac{1}{2}
\left[B_\beta(\lambda)^{-1/2}(D_v\otimes1)B_\beta(\lambda)^{-1/2}\right]
B_\beta(\lambda)^\frac{1}{2}.
\end{equation*}
Since $\Lambda_\beta$ is convex, \eqref{eq:limiting-hessian} proves strict convexity there.

Suppose that $\lambda,\mu\in\overline\Lambda_\beta$ have finite values $\widehat\phi_\beta(\lambda)=\widehat\phi_\beta(\mu)=F(\beta)$. For $0<\delta\leq1$, the segment joining $\lambda+\delta\1$ and $\mu+\delta\1$ lies in $\Lambda_\beta$, and the norm of $B_\beta$ on this segment is at most $C:=1+\max\{\|\lambda\|_\infty,\|\mu\|_\infty\}+2\beta\|\mathsf M\|$. Integrating the Hessian bound gives
\begin{equation}\label{eq:regularized-midpoint}
F(\beta)\leq\phi_\beta\left(\frac{\lambda+\mu}{2}+\delta\1\right)
\stackrel{\eqref{eq:limiting-hessian}}{\leq}
\frac{\phi_\beta(\lambda+\delta\1)+\phi_\beta(\mu+\delta\1)}2
-\frac{\|\lambda-\mu\|^2}{16C^2}.
\end{equation}
Monotone convergence in the scalar spectral integrals defining \eqref{eq:two-sided-objective} gives
\begin{equation*}
\lim_{\delta\downarrow0}\phi_\beta(\nu+\delta\1)
=\widehat\phi_\beta(\nu)=F(\beta)
\qquad\text{for }\nu\in\{\lambda,\mu\}.
\end{equation*}
Letting $\delta\downarrow0$ in \eqref{eq:regularized-midpoint} yields $\|\lambda-\mu\|^2/(16C^2)\leq0$, so $\lambda=\mu$.
\end{proof}

Fix $\beta>0$ and suppose that $\lambda_\beta\in\Lambda_\beta$ satisfies $\phi_\beta(\lambda_\beta)=F(\beta)$. By Lemma~\ref{lem:interior-minimizer-uniqueness}, this $\lambda_\beta$ is unique. Define the positive definite matrix $V_\beta\in M_d(\R)$ by
\begin{equation}\label{eq:block-norm-covariance}
(V_\beta)_{ab}:=
2d\tau_d
\left[(E_{aa}\otimes1)B_\beta(\lambda_\beta)^{-1}(E_{bb}\otimes1)B_\beta(\lambda_\beta)^{-1}\right].
\end{equation}
Thus $V_\beta=4\nabla^2\phi_\beta(\lambda_\beta)$.

Choose a compact convex neighborhood $U$ of $\lambda_\beta$ with $U\subset\Lambda_\beta$, constants $0<\varepsilon<T<\infty$ such that
\begin{equation*}
\operatorname{spec}(D_\lambda\otimes1-2\beta\mathsf M)\subset(\varepsilon,T)
\quad\text{for every }\lambda\in U,
\end{equation*}
and a real-valued $f\in C_c^\infty(\R)$ that agrees with $\log$ on a neighborhood of $[\varepsilon,T]$. Define
\begin{equation}\label{eq:truncated-interior-objective}
\widetilde\phi_{N,\beta}(\lambda):=
\frac12\sum_{a=1}^d(\lambda_a-1)
-\frac1{2N}\Tr f(B_{N,\beta}(\lambda)).
\end{equation}
Choose a deterministic minimizer $\bar\lambda_{N,\beta}\in\operatorname*{argmin}_{\lambda\in U}\E\widetilde\phi_{N,\beta}(\lambda)$, which is unique for all sufficiently large $N$ by Lemma~\ref{lem:coarea-minimizer-expansion}, and set
\begin{align}
X_{N,\beta}&:=\Tr f(B_{N,\beta}(\bar\lambda_{N,\beta}))
-\E\Tr f(B_{N,\beta}(\bar\lambda_{N,\beta})),\label{eq:interior-centered-statistic}\\
m_{N,\beta}&:=N\left(\E\widetilde\phi_{N,\beta}(\bar\lambda_{N,\beta})-F(\beta)\right)
+\frac d2\log2-\frac12\log\det V_\beta.\label{eq:interior-centering}
\end{align}
The additional assumption for the next theorem is
\begin{equation}\label{eq:interior-resolvent-variance}
\sup_{N\in\mathbb N}\sup_{\lambda\in U}\max_{a\in\{1,\ldots,d\}}
\operatorname{Var}\!\left(\Tr[P_af'(B_{N,\beta}(\lambda))]\right)<\infty.
\end{equation}

\begin{theorem}\label{thm:interior-fluctuations}
Assume \eqref{eq:interior-resolvent-variance}. Then $m_{N,\beta}=o(N)$ and
\begin{equation}\label{eq:interior-logdet-reduction}
N(F_N(\beta)-F(\beta))-m_{N,\beta}
=-\frac12X_{N,\beta}+o_{\mathbb P}(1).
\end{equation}
Consequently, if $X_{N,\beta}\xrightarrow[N\to\infty]{\mathrm d}X_\beta$, then
\begin{equation*}
N(F_N(\beta)-F(\beta))-m_{N,\beta}
\xrightarrow[N\to\infty]{\mathrm d}-\frac12X_\beta.
\end{equation*}
\end{theorem}

\begin{lemma}\label{lem:coarea-minimizer-expansion}
Let $\lambda_{N,\beta}$ be the unique minimizer of $\phi_{N,\beta}$ over the domain where $B_{N,\beta}(\lambda)\succ0$. For the fixed $\lambda_\beta\in\Lambda_\beta$ with $\phi_\beta(\lambda_\beta)=F(\beta)$, one has
\begin{equation}\label{eq:mean-objective-c2}
\lim_{N\to\infty}\|\E\widetilde\phi_{N,\beta}-\phi_\beta\|_{C^2(U)}=0.
\end{equation}
For all sufficiently large $N$, the minimizer $\bar\lambda_{N,\beta}$ is unique, and
\begin{equation}\label{eq:deterministic-minimizer-convergence}
\lim_{N\to\infty}\|\bar\lambda_{N,\beta}-\lambda_\beta\|=0.
\end{equation}
Moreover,
\begin{equation}\label{eq:finite-minimizer-convergence}
\lim_{N\to\infty}\|\lambda_{N,\beta}-\lambda_\beta\|=0,
\end{equation}
in probability and
\begin{equation}\label{eq:exact-spherical-expansion}
\log\int_{\mathcal S_N}e^{\beta H_N(\sigma)}\,\d\mu_N^{\otimes d}(\sigma)
=N\phi_{N,\beta}(\lambda_{N,\beta})
+\frac d2\log2-\frac12\log\det V_\beta+o_{\mathbb P}(1).
\end{equation}
\end{lemma}

\begin{proof}
\smallskip
\noindent\emph{Step 1: Convergence of the derivatives and minimizers.}
We use the norm
\begin{equation*}
\|u\|_{C^2(U)}:=\max_{|\alpha|\leq2}\sup_{\lambda\in U}|\partial^\alpha u(\lambda)|,
\end{equation*}
where $\alpha$ ranges over multi-indices of order at most two. Work on the probability-one event in Proposition~\ref{prop:strong-convergence}. Choose a compact neighborhood $U'$ of $U$ contained in $\Lambda_\beta$. Equations \eqref{eq:affine-spectrum}--\eqref{eq:affine-trace} give uniform convergence of $\phi_{N,\beta}$ to $\phi_\beta$ on $U'$, and a uniform positive lower bound for $B_{N,\beta}(\lambda)$ on $U'$. Differentiating the inverse in the logarithmic determinant therefore bounds all derivatives through order three uniformly in $N$ and $\lambda\in U'$, for all sufficiently large $N$. For $1\leq|\alpha|\leq2$, finite differences on $U'$ give, for all sufficiently small $h>0$,
\begin{equation*}
\sup_{\lambda\in U}|\partial^\alpha(\phi_{N,\beta}-\phi_\beta)(\lambda)|
\leq C h^{-|\alpha|}\sup_{\lambda\in U'}|\phi_{N,\beta}(\lambda)-\phi_\beta(\lambda)|+Ch.
\end{equation*}
Taking $N\to\infty$ and then $h\downarrow0$ proves
\begin{equation}\label{eq:interior-objective-c2}
\lim_{N\to\infty}\|\phi_{N,\beta}-\phi_\beta\|_{C^2(U)}=0
\quad\text{almost surely}.
\end{equation}
The choice of $f$ and \eqref{eq:affine-spectrum} give $\widetilde\phi_{N,\beta}=\phi_{N,\beta}$ on $U$ for all sufficiently large $N$. The divided-difference formula for the Fr\'echet derivative of $f'$ gives
\begin{equation*}
\|Df'(B)[A]\|_{\mathrm{HS}}\leq\|f''\|_\infty\|A\|_{\mathrm{HS}}
\quad\text{for every }B=B^*,\ A=A^*.
\end{equation*}
Since $\partial_a\Tr f(B_{N,\beta}(\lambda))=\Tr[P_af'(B_{N,\beta}(\lambda))]$, this bound and $\|P_a\|_{\mathrm{HS}}^2=N$ yield the deterministic estimate
\begin{equation*}
\sup_{N\in\mathbb N}\|\widetilde\phi_{N,\beta}\|_{C^2(U)}\leq C_{U,f,d}.
\end{equation*}
Dominated convergence now proves \eqref{eq:mean-objective-c2}, including differentiation under the expectation. Equation \eqref{eq:limiting-hessian} and compactness of $U$ give $\inf_{\lambda\in U}\lambda_{\min}(\nabla^2\phi_\beta(\lambda))>0$. Thus \eqref{eq:mean-objective-c2} makes $\E\widetilde\phi_{N,\beta}$ strictly convex on $U$ for all sufficiently large $N$. Its uniform convergence to $\phi_\beta$ and Lemma~\ref{lem:interior-minimizer-uniqueness} prove uniqueness of $\bar\lambda_{N,\beta}$ and \eqref{eq:deterministic-minimizer-convergence}.

The stationarity equations for \eqref{eq:finite-precision-objective} are
\begin{equation}\label{eq:fluctuation-stationarity}
\Tr\left[P_aB_{N,\beta}(\lambda_{N,\beta})^{-1}\right]=N
\quad\text{for every }a\in\{1,\ldots,d\}.
\end{equation}
Equation \eqref{eq:lambda-stationarity} bounds $\lambda_{N,\beta}$, and $B_{N,\beta}(\lambda_\beta)\succ0$ for all sufficiently large $N$. On the segment joining $\lambda_\beta$ and $\lambda_{N,\beta}$, the precision matrix is positive definite and has norm bounded by some $C<\infty$, independent of $N$. The finite version of \eqref{eq:limiting-hessian} gives
\begin{equation*}
v^{\mathsf T}\nabla^2\phi_{N,\beta}(\lambda)v
=\frac1{2N}\left\|B_{N,\beta}(\lambda)^{-1/2}(D_v\otimes I_N)B_{N,\beta}(\lambda)^{-1/2}\right\|_{\mathrm{HS}}^2
\geq\frac{\|v\|^2}{2C^2}.
\end{equation*}
Equation \eqref{eq:interior-objective-c2} gives $\lim_{N\to\infty}\|\nabla\phi_{N,\beta}(\lambda_\beta)\|=\|\nabla\phi_\beta(\lambda_\beta)\|=0$. Integrating the Hessian along the segment and using \eqref{eq:fluctuation-stationarity}, we obtain
\begin{equation*}
\|\lambda_{N,\beta}-\lambda_\beta\|
\leq2C^2\|\nabla\phi_{N,\beta}(\lambda_\beta)\|.
\end{equation*}
This proves \eqref{eq:finite-minimizer-convergence}, including exclusion of a boundary subsequential limit. Consequently, with probability tending to one as $N\to\infty$, both $B_{N,\beta}(\lambda_{N,\beta})$ and its inverse are bounded by deterministic positive scalars.

\smallskip
\noindent\emph{Step 2: An exact coarea identity.}
Condition on $M_N$, set $\Sigma_N:=B_{N,\beta}(\lambda_{N,\beta})^{-1}$, and let $X$ be a centered Gaussian vector in $\R^{dN}$ with covariance $\Sigma_N$. Put
\begin{equation*}
Y_a:=\|P_aX\|^2,
\end{equation*}
and let $p_N$ denote the joint density of $(Y_1,\ldots,Y_d)$. By \eqref{eq:fluctuation-stationarity}, we have $\E_XY_a=N$. The coarea formula applied to
\begin{equation*}
x\longmapsto(\|P_1x\|^2,\ldots,\|P_dx\|^2)
\end{equation*}
gives
\begin{equation}\label{eq:coarea-identity}
\begin{aligned}
\int_{\mathcal S_N}e^{\beta H_N(\sigma)}\,\d\mu_N^{\otimes d}(\sigma)
&=\exp\left(\frac N2\sum_{a=1}^d(\lambda_{N,\beta})_a\right)
\det B_{N,\beta}(\lambda_{N,\beta})^{-1/2}\\
&\qquad\times\left[2^{N/2}\Gamma(N/2)N^{1-N/2}\right]^d p_N(N\1).
\end{aligned}
\end{equation}
Indeed, the Jacobian of $x\mapsto(\|P_1x\|^2,\ldots,\|P_dx\|^2)$ on the level set $(N,\ldots,N)$ is $2^dN^{d/2}$, while the surface area of the sphere of radius $\sqrt N$ in $\R^N$ is
\begin{equation*}
\frac{2\pi^{N/2}N^{(N-1)/2}}{\Gamma(N/2)}.
\end{equation*}
These constants give \eqref{eq:coarea-identity}.

\smallskip
\noindent\emph{Step 3: A local central limit theorem for the block norms.}
In this proof, set
\begin{equation*}
(V_N)_{ab}:=\frac2N\Tr\left[P_aB_{N,\beta}(\lambda_{N,\beta})^{-1}P_bB_{N,\beta}(\lambda_{N,\beta})^{-1}\right].
\end{equation*}
Equations \eqref{eq:finite-minimizer-convergence} and \eqref{eq:interior-objective-c2}, together with $V_N=4\nabla^2\phi_{N,\beta}(\lambda_{N,\beta})$ and \eqref{eq:block-norm-covariance}, give
\begin{equation}\label{eq:VN-limit}
\lim_{N\to\infty}\|V_N-V_\beta\|=0,
\end{equation}
in probability. For $t\in\R^d$, the conditional characteristic function of $(Y-N\1)/\sqrt N$ is
\begin{equation}\label{eq:block-norm-characteristic}
\E_X\exp\left(\frac{\mathrm i}{\sqrt N}\sum_{a=1}^dt_a(Y_a-N)\right)
=e^{-\mathrm i\sqrt N\sum_at_a}
\det\left(I_{dN}-\frac{2\mathrm i}{\sqrt N}B_{N,\beta}(\lambda_{N,\beta})^{-1}(D_t\otimes I_N)\right)^{-1/2}.
\end{equation}
Expanding its logarithm, the linear term cancels by \eqref{eq:fluctuation-stationarity}, and, uniformly for $t$ in compact sets,
\begin{equation}\label{eq:block-norm-characteristic-limit}
\log\E_X\exp\left(\frac{\mathrm i}{\sqrt N}\sum_{a=1}^dt_a(Y_a-N)\right)
=-\frac12t^{\mathsf T}V_Nt+O(N^{-1/2}|t|^3).
\end{equation}
On the event that the eigenvalues of the covariance are bounded above and below by deterministic positive constants, the characteristic functions admit an integrable tail bound. Indeed, if $A_N(t):=\Sigma_N^\frac{1}{2}(D_t\otimes I_N)\Sigma_N^\frac{1}{2}$, then the absolute value of the right-hand side of \eqref{eq:block-norm-characteristic} is
\begin{equation*}
\det\left(I_{dN}+\frac4N A_N(t)^2\right)^{-1/4}.
\end{equation*}
If $cI_{dN}\preceq\Sigma_N\preceq CI_{dN}$, the singular-value inequality for multiplication by the invertible matrix $\Sigma_N^\frac{1}{2}$ gives
\begin{equation*}
s_j(A_N(t))\geq c\,s_j(D_t\otimes I_N)
\quad\text{for every }j\in\{1,\ldots,dN\}.
\end{equation*}
Since the singular values of $D_t\otimes I_N$ are $|t_a|$, each repeated $N$ times, we obtain
\begin{equation}\label{eq:block-characteristic-tail}
\left|\E_X\exp\left(\frac{\mathrm i}{\sqrt N}\sum_{a=1}^dt_a(Y_a-N)\right)\right|
\leq\prod_{a=1}^d\left(1+\frac{c t_a^2}{N}\right)^{-N/4}.
\end{equation}
Suppose that $cI_{dN}\preceq\Sigma_N\preceq CI_{dN}$. For some $c_0>0$, if $|t|\leq c_0\sqrt N$, then $\log(1+x)\geq c x$ on the relevant bounded interval and \eqref{eq:block-characteristic-tail} is bounded by $e^{-c|t|^2}$. Hence the integral over $R<|t|\leq c_0\sqrt N$ tends to zero uniformly in $N$ as $R\to\infty$. For $|t|>c_0\sqrt N$, some coordinate satisfies $|t_a|>c_1\sqrt N$. For some $N_0$, we have
\begin{equation*}
\sup_{N\geq N_0}\int_{\R}\left(1+\frac{c x^2}{N}\right)^{-N/4}\,\d x<\infty,
\end{equation*}
while, for all sufficiently large $N$,
\begin{equation*}
\int_{|x|>c_1\sqrt N}\left(1+\frac{c x^2}{N}\right)^{-N/4}\,\d x
=\sqrt N\int_{|u|>c_1}(1+cu^2)^{-N/4}\,\d u
\leq C\sqrt N e^{-c'N}.
\end{equation*}
Therefore, splitting according to the largest coordinate in \eqref{eq:block-characteristic-tail}, the integral over $|t|>c_0\sqrt N$ also tends to zero. Since the event $cI_{dN}\preceq\Sigma_N\preceq CI_{dN}$ has probability tending to one as $N\to\infty$, we obtain, for every $\delta>0$,
\begin{equation}\label{eq:block-characteristic-uniform-tail}
\lim_{R\to\infty}\limsup_{N\to\infty}
\P\left(
\int_{|t|>R}
\left|\E_X\exp\left(\frac{\mathrm i}{\sqrt N}\sum_{a=1}^dt_a(Y_a-N)\right)\right|\,\mathrm dt
>\delta
\right)=0.
\end{equation}
Fourier inversion, \eqref{eq:VN-limit}, \eqref{eq:block-norm-characteristic-limit}, and \eqref{eq:block-characteristic-uniform-tail} therefore give
\begin{equation}\label{eq:block-norm-local-clt}
p_N(N\1)
=N^{-d/2}(2\pi)^{-d/2}(\det V_\beta)^{-1/2}
\bigl(1+o_{\mathbb P}(1)\bigr).
\end{equation}

\smallskip
\noindent\emph{Step 4: Proof of \eqref{eq:exact-spherical-expansion}.}
Stirling's formula gives
\begin{equation}\label{eq:stirling-spherical-correction}
\frac{dN}{2}+d\log\Gamma(N/2)+\frac{dN}{2}\log2+d(1-N/2)\log N
=\frac d2\log(4\pi N)+o(1).
\end{equation}
Combining \eqref{eq:coarea-identity}, \eqref{eq:block-norm-local-clt}, and \eqref{eq:stirling-spherical-correction} proves \eqref{eq:exact-spherical-expansion}.
\end{proof}

\begin{proof}[Proof of Theorem~\ref{thm:interior-fluctuations}]
Recall $\widetilde\phi_{N,\beta}$ from \eqref{eq:truncated-interior-objective}. By Lemma~\ref{lem:coarea-minimizer-expansion}, the deterministic point $\bar\lambda_{N,\beta}$ lies in the interior of $U$ for all sufficiently large $N$, and $\nabla\E\widetilde\phi_{N,\beta}(\bar\lambda_{N,\beta})=0$. Hence, for every $a\in\{1,\ldots,d\}$,
\begin{align*}
\partial_a\widetilde\phi_{N,\beta}(\bar\lambda_{N,\beta})
=-\frac1{2N}\left\{
\Tr[P_af'(B_{N,\beta}(\bar\lambda_{N,\beta}))]
-\E\Tr[P_af'(B_{N,\beta}(\bar\lambda_{N,\beta}))]\right\}\stackrel{\eqref{eq:interior-resolvent-variance}}{=}O_{\mathbb P}(N^{-1}).
\end{align*}
By \eqref{eq:affine-spectrum} and \eqref{eq:finite-minimizer-convergence}, with probability tending to one as $N\to\infty$, $\lambda_{N,\beta}\in U$ and $\varepsilon I_{dN}\preceq B_{N,\beta}(\lambda)\preceq TI_{dN}$ for every $\lambda\in U$. On this event, $\widetilde\phi_{N,\beta}=\phi_{N,\beta}$ on $U$, and the finite version of \eqref{eq:limiting-hessian} gives
\begin{equation}\label{eq:interior-hessian-bounds}
\frac{\|v\|^2}{2T^2}
\leq v^{\mathsf T}\nabla^2\phi_{N,\beta}(\lambda)v
\leq\frac{\|v\|^2}{2\varepsilon^2}\quad \text{for every $\lambda\in U,\ v\in\R^d$}.
\end{equation}
Integrating the lower bound along the segment joining the two minimizers yields
\begin{equation}\label{eq:finite-deterministic-minimizer-distance}
\|\lambda_{N,\beta}-\bar\lambda_{N,\beta}\|
\leq2T^2\|\nabla\widetilde\phi_{N,\beta}(\bar\lambda_{N,\beta})\|
=O_{\mathbb P}(N^{-1}).
\end{equation}
Taylor's theorem at $\lambda_{N,\beta}$, where the gradient vanishes, together with \eqref{eq:finite-deterministic-minimizer-distance} and the upper bound in \eqref{eq:interior-hessian-bounds}, gives
\begin{equation}\label{eq:minimizer-value-remainder}
N\left|\phi_{N,\beta}(\lambda_{N,\beta})
-\widetilde\phi_{N,\beta}(\bar\lambda_{N,\beta})\right|
\leq\frac{N}{4\varepsilon^2}\|\lambda_{N,\beta}-\bar\lambda_{N,\beta}\|^2
=O_{\mathbb P}(N^{-1}).
\end{equation}
Combining \eqref{eq:exact-spherical-expansion}, \eqref{eq:minimizer-value-remainder}, and the definitions \eqref{eq:interior-centered-statistic}--\eqref{eq:interior-centering} proves \eqref{eq:interior-logdet-reduction}. Equations \eqref{eq:mean-objective-c2} and \eqref{eq:deterministic-minimizer-convergence} give $m_{N,\beta}=o(N)$. The distributional assertion follows from Slutsky's theorem.
\end{proof}

We now assume that $g_{11}$ is Gaussian and use \eqref{eq:gaussian-block-pencil}. Recall the function $f$ and the matrices $B_{N,\beta}(\lambda)$ defined before Theorem~\ref{thm:interior-fluctuations}. A polynomial trace of $B_{N,\beta}(\lambda_\beta)$ expands, through its $d\times d$ block indices, into a finite linear combination of traces of words in the independent real Ginibre matrices $G_{r,N}$ and their transposes, with deterministic coefficients depending only on the matrices $A_r$. The Wick expansion and real second-order limit distribution in \cite[Sections 4--6]{redelmeier} therefore give a second-order limit distribution for $B_{N,\beta}(\lambda_\beta)$.

For every compact $U\subset\Lambda_\beta$, there are deterministic constants $0<\varepsilon_U<T_U<\infty$ and $c_U,C_U>0$ such that
\begin{equation}\label{eq:gaussian-uniform-spectrum}
\P\left(\operatorname{spec}B_{N,\beta}(\lambda)\not\subset[\varepsilon_U,T_U]
\text{ for some }\lambda\in U\right)\leq C_Ue^{-c_UN}.
\end{equation}
Indeed, each extreme eigenvalue is $O(N^{-1/2})$-Lipschitz in the standard Gaussian entries. Strong convergence and \eqref{eq:ginibre-norm-moments} give convergence of its expectation at a fixed $\lambda$. Gaussian concentration and a finite net of $U$, using $\|B_{N,\beta}(\lambda)-B_{N,\beta}(\mu)\|=\|\lambda-\mu\|_\infty$, prove \eqref{eq:gaussian-uniform-spectrum}. For the particular $U$, $\varepsilon$, and $T$ chosen before Theorem~\ref{thm:interior-fluctuations}, the concentration-and-net proof of \eqref{eq:gaussian-uniform-spectrum} applies with $\varepsilon_U=\varepsilon$ and $T_U=T$.

In particular, \eqref{eq:gaussian-uniform-spectrum} gives $\lim_{N\to\infty}(dN)^8\P(\|B_{N,\beta}(\lambda_\beta)\|>T)=0$, which verifies Assumption A1 of \cite[Section 3]{diaz-mingo}. Assumption A2 follows from the Gaussian Poincar\'e inequality. Set
\begin{equation*}
H_{ijr}:=\frac{\partial}{\partial z_{ij}^r}B_{N,\beta}(\lambda)
=-\frac{\beta}{\sqrt N}\left(A_r\otimes E_{ij}+A_r^{\mathsf T}\otimes E_{ji}\right).
\end{equation*}
Since $d$, $m$, and the matrices $A_r$ are fixed, summing the squared matrix entries gives, for every matrix $K$ of size $dN$,
\begin{equation}\label{eq:gaussian-hs-estimate}
\sum_{r=1}^m\sum_{i,j=1}^N|\Tr(KH_{ijr})|^2
\leq\frac{C}{N}\|K\|_{\mathrm{HS}}^2.
\end{equation}
For any $u\in C^1(\R)$ with bounded derivative, $\partial_{z_{ij}^r}\Tr u(B_{N,\beta}(\lambda))=\Tr[u'(B_{N,\beta}(\lambda))H_{ijr}]$. Hence \eqref{eq:gaussian-hs-estimate}, applied with $K=u'(B_{N,\beta}(\lambda))$, and $\|u'(B_{N,\beta}(\lambda))\|_{\mathrm{HS}}^2\leq dN\|u'\|_\infty^2$ give, for a constant $C$ independent of $N$ and $u$,
\begin{equation}\label{eq:logdet-poincare-bound}
\sum_{r=1}^m\sum_{i,j=1}^N
\E\left|\frac{\partial}{\partial z_{ij}^r}
\Tr u(B_{N,\beta}(\lambda_\beta))\right|^2
\leq C\|u'\|_\infty^2.
\end{equation}
Consequently, \cite[Theorem 6(b)]{diaz-mingo} gives the finite limit
\begin{equation}\label{eq:logdet-variance}
v_\beta:=\lim_{N\to\infty}
\operatorname{Var}\left(\Tr f(B_{N,\beta}(\lambda_\beta))\right),
\end{equation}
and \cite[Proposition 12]{diaz-mingo} gives
\begin{equation*}
\Tr f(B_{N,\beta}(\lambda_\beta))
-\E\Tr f(B_{N,\beta}(\lambda_\beta))
\xrightarrow[N\to\infty]{\mathrm d}\mathcal N(0,v_\beta).
\end{equation*}
For two choices of $f$ agreeing with $\log$ near the limiting spectrum, apply \eqref{eq:gaussian-uniform-spectrum} on a smaller common spectral interval. Their traces differ by at most $O(N)$ on an event of probability $O(e^{-cN})$, and hence their difference tends to zero in $L^2$. Thus \eqref{eq:logdet-variance} is independent of the choice of $f$.

To apply the contour formula in \cite[Theorem 11]{diaz-mingo}, set
\begin{equation*}
a_0:=\frac{T+\varepsilon}{2}
\qquad\text{and}\qquad
\widehat B_{N,\beta}:=B_{N,\beta}(\lambda_\beta)-a_0I_{dN}.
\end{equation*}
On the event $\operatorname{spec}B_{N,\beta}(\lambda_\beta)\subset[\varepsilon,T]$, we have $\|\widehat B_{N,\beta}\|\leq(T-\varepsilon)/2<T/2<a_0$. Therefore, \eqref{eq:gaussian-uniform-spectrum} verifies Assumption A1 for $\widehat B_{N,\beta}$ with the bound $T/2$. Choose a real $\widehat f\in C_c^\infty(\R)$ that agrees with $\log(a_0+x)$ on a neighborhood of $[-T/2,T/2]$. Its restriction to this interval extends analytically as $\log(a_0+z)$ to $\{z\in\C:\operatorname{Re}z>-a_0\}$. The deterministic scalar shift preserves the second-order polynomial limit and the Gaussian Poincar\'e bound, so \cite[Theorem 11]{diaz-mingo} applies to $(\widehat B_{N,\beta},\widehat f)$. Moreover, the two traces agree on the spectral event, and
\begin{equation*}
\E\left|\Tr \widehat f(\widehat B_{N,\beta})-\Tr f(B_{N,\beta}(\lambda_\beta))\right|^2
\leq(dN)^2(\|\widehat f\|_\infty+\|f\|_\infty)^2C_Ue^{-c_UN}.
\end{equation*}
Since $\lim_{N\to\infty}N^2e^{-c_UN}=0$, we obtain $\lim_{N\to\infty}\operatorname{Var}(\Tr \widehat f(\widehat B_{N,\beta}))=v_\beta$ by \eqref{eq:logdet-poincare-bound} and \eqref{eq:logdet-variance}. Thus the shifted contour formula determines the same $v_\beta$. For complex Gaussian block matrices, \cite[Theorem 2]{diaz-mingo-belinschi} gives the analogous second-order Cauchy-transform formula.

\begin{corollary}\label{cor:gaussian-interior-fluctuations}
If $g_{11}$ is Gaussian, then the variance bound \eqref{eq:interior-resolvent-variance} holds, $m_{N,\beta}=o(N)$, and
\begin{equation}\label{eq:interior-centered-clt}
N(F_N(\beta)-F(\beta))-m_{N,\beta}
\xrightarrow[N\to\infty]{\mathrm d}
\mathcal N\left(0,\frac{v_\beta}{4}\right).
\end{equation}
If $\lim_{N\to\infty}m_{N,\beta}=m_\beta$, then
\begin{equation*}
N(F_N(\beta)-F(\beta))
\xrightarrow[N\to\infty]{\mathrm d}
\mathcal N\left(m_\beta,\frac{v_\beta}{4}\right).
\end{equation*}
\end{corollary}

\begin{proof}
The Gaussian Poincar\'e inequality also applies to weighted traces. By self-adjointness of the Fr\'echet derivative with respect to the Hilbert--Schmidt pairing,
\begin{equation*}
\frac{\partial}{\partial z_{ij}^r}\Tr[P_af'(B_{N,\beta}(\lambda))]
=\Tr\left[Df'(B_{N,\beta}(\lambda))[P_a]H_{ijr}\right].
\end{equation*}
Applying \eqref{eq:gaussian-hs-estimate} with $K=Df'(B_{N,\beta}(\lambda))[P_a]$ and the divided-difference bound $\|Df'(B)[P_a]\|_{\mathrm{HS}}\leq\|f''\|_\infty\|P_a\|_{\mathrm{HS}}$, gives, for every $a\in\{1,\ldots,d\}$,
\begin{equation*}
\sup_{N\in\mathbb N}\sup_{\lambda\in U}
\operatorname{Var}\!\left(\Tr[P_af'(B_{N,\beta}(\lambda))]\right)
\leq C\|f''\|_\infty^2\frac{\|P_a\|_{\mathrm{HS}}^2}{N}
=C\|f''\|_\infty^2,
\end{equation*}
which verifies \eqref{eq:interior-resolvent-variance}.

The parameter $\bar\lambda_{N,\beta}$ is deterministic and satisfies \eqref{eq:deterministic-minimizer-convergence}. For each fixed polynomial trace statistic, expanding $B_{N,\beta}(\lambda)=D_\lambda\otimes I_N-2\beta M_N$ and applying Wick's formula expresses its expectations and covariances as finite sums over Gaussian pairings; each summand is a polynomial in the entries of $D_\lambda$ and the fixed matrices $A_1,\ldots,A_m$. Since \eqref{eq:deterministic-minimizer-convergence} gives $\lim_{N\to\infty}\|\bar\lambda_{N,\beta}-\lambda_\beta\|=0$, these first- and second-order limits are unchanged when $\lambda_\beta$ is replaced by $\bar\lambda_{N,\beta}$. Hence $B_{N,\beta}(\bar\lambda_{N,\beta})$ has the same second-order limit distribution as $B_{N,\beta}(\lambda_\beta)$. The spectral estimate \eqref{eq:gaussian-uniform-spectrum} and the derivative estimate \eqref{eq:logdet-poincare-bound} are uniform when $\lambda_\beta$ is replaced by $\lambda\in U$. Applying \cite[Theorem 6(b) and Proposition 12]{diaz-mingo} to $B_{N,\beta}(\bar\lambda_{N,\beta})$ therefore gives
\begin{equation*}
X_{N,\beta}\xrightarrow[N\to\infty]{\mathrm d}\mathcal N(0,v_\beta).
\end{equation*}
Theorem~\ref{thm:interior-fluctuations} proves \eqref{eq:interior-centered-clt}. The final assertion follows from Slutsky's theorem.
\end{proof}

\subsection{Transfer of GOE Tracy--Widom fluctuations when \texorpdfstring{$\lambda_\beta\in\partial\Lambda_\beta$ and $\theta_\beta>0$}{lambda beta is on the boundary and theta beta is positive}}\label{sec:boundary-fluctuations}

\begin{proof}[Proof of Proposition~\ref{prop:boundary-resolvent-multiplier}]
Recall $g_\beta$ and $\widehat\phi_\beta$ from \eqref{eq:edge-function} and \eqref{eq:two-sided-objective}. Differentiating \eqref{eq:edge-shift} gives
\begin{equation}\label{eq:edge-normal-shift}
\nabla g_\beta(\lambda_\beta)\cdot\1=-1.
\end{equation}
For $v\in\R^d$ with $\nabla g_\beta(\lambda_\beta)\cdot v<0$, one has $\lambda_\beta+sv\in\Lambda_\beta$ for all sufficiently small $s>0$. The minimizing property in Definition~\ref{def:boundary-multiplier} then gives
\begin{equation*}
\nabla g_\beta(\lambda_\beta)\cdot v<0
\quad\Longrightarrow\quad
\nabla\widehat\phi_\beta(\lambda_\beta)\cdot v\geq0.
\end{equation*}
If $\nabla g_\beta(\lambda_\beta)\cdot v=0$, apply this implication to $v+\delta\1$ and $-v+\delta\1$ and let $\delta\downarrow0$. Thus $\nabla\widehat\phi_\beta(\lambda_\beta)$ vanishes on the kernel of $v\mapsto\nabla g_\beta(\lambda_\beta)\cdot v$. Taking $v=\1$ and using \eqref{eq:edge-normal-shift} proves that \eqref{eq:boundary-multiplier} has the unique solution $\theta_\beta=\nabla\widehat\phi_\beta(\lambda_\beta)\cdot\1\geq0$.

Let $\nu_0$ be the spectral distribution of $B_\beta(\lambda_\beta)=D_{\lambda_\beta}\otimes1-2\beta\mathsf M\succeq0$ under $\tau_d$. Finiteness of $\widehat\phi_\beta(\lambda_\beta)$ implies $\nu_0(\{0\})=0$ and integrability of $\log x$ near zero. For $t>0$ and $x>0$, put $q_t(x):=(\log(x+t)-\log x)/t$. Since $q_t(x)=\int_0^1(x+st)^{-1}\,\d s$, one has $q_{t'}(x)\geq q_t(x)$ for $0<t'<t$ and $\lim_{t\downarrow0}q_t(x)=x^{-1}$. Consequently, differentiability and monotone convergence give
\begin{equation}\label{eq:boundary-inverse-trace}
\begin{aligned}
\theta_\beta
&=\lim_{t\downarrow0}\frac{\widehat\phi_\beta(\lambda_\beta+t\1)-\widehat\phi_\beta(\lambda_\beta)}{t}\\
&=\lim_{t\downarrow0}\left(\frac d2-\frac d2\int q_t(x)\,\nu_0(\d x)\right)
=\frac d2-\frac d2\int_{(0,\infty)}\frac{\nu_0(\d x)}{x}.
\end{aligned}
\end{equation}
Equation \eqref{eq:boundary-inverse-trace} implies $\int_{(0,\infty)}x^{-1}\,\nu_0(\d x)<\infty$. The matrices $Q_{\lambda_\beta+t\1}$ are positive definite and increase in the positive semidefinite order as $t\downarrow0$. Moreover,
\begin{equation*}
\Tr Q_{\lambda_\beta+t\1}
=d\int_{(0,\infty)}\frac{\nu_0(\d x)}{x+t}
\leq d\int_{(0,\infty)}\frac{\nu_0(\d x)}{x}<\infty.
\end{equation*}
Monotonicity and this uniform trace bound imply convergence in matrix norm to the finite matrix $Q^\partial_{\lambda_\beta}$ defined in \eqref{eq:boundary-resolvent}. Taking $t\downarrow0$ identifies its trace as $d\int x^{-1}\,\nu_0(\d x)$. Equation \eqref{eq:boundary-inverse-trace} proves \eqref{eq:boundary-multiplier-trace}.
\end{proof}

Fix $\beta>0$ and $\lambda_\beta\in\partial\Lambda_\beta$ satisfying Definition~\ref{def:boundary-multiplier}, and suppose that $\theta_\beta>0$. The next definition specifies what it means for $(g_\beta,\widehat\phi_\beta)$ to be nondegenerate at $\lambda_\beta$.

\begin{definition}[Nondegeneracy at $\lambda_\beta\in\partial\Lambda_\beta$]\label{def:boundary-nondegeneracy}
We say that $(g_\beta,\widehat\phi_\beta)$ is \textit{nondegenerate at $\lambda_\beta$} if $g_\beta$ is twice continuously differentiable near $\lambda_\beta$ and there are a neighborhood $U$ of $\lambda_\beta$ on which $\widehat\phi_\beta$ is finite, neighborhoods $W\subset\R^{d-1}$ of $0$ and $I\subset\R$ of $0$, a twice continuously differentiable diffeomorphism $\Gamma:W\times I\to U$, and a positive definite matrix $H_\beta\in M_{d-1}(\R)$ such that
\begin{equation*}
\Gamma(0,0)=\lambda_\beta
\qquad\text{and}\qquad
g_\beta(\Gamma(y,r))=r,
\end{equation*}
and, as $(y,r)\to(0,0)$,
\begin{equation}\label{eq:boundary-local-expansion}
\widehat\phi_\beta(\Gamma(y,r))
=F(\beta)-\theta_\beta r
+\frac12 y^{\mathsf T}H_\beta y
+o\bigl(|r|+\|y\|^2\bigr).
\end{equation}
Since $\widehat\phi_\beta=\phi_\beta$ on $\Lambda_\beta$, \eqref{eq:boundary-local-expansion} is consistent with \eqref{eq:boundary-multiplier}. It also implies that, after decreasing $U$ if necessary,
\begin{equation*}
\widehat\phi_\beta(\lambda)>\widehat\phi_\beta(\lambda_\beta)=F(\beta)
\end{equation*}
for every $\lambda\in(U\cap\overline{\Lambda}_\beta)\setminus\{\lambda_\beta\}$. When $d=1$, the variables $y$ and the matrix $H_\beta$ are absent.
\end{definition}

\begin{remark}
When the spectral density vanishes at square-root rate at an edge, the two-sided logarithmic potential in \eqref{eq:two-sided-objective} is typically continuously differentiable in the normal direction, while its second normal derivative is singular. Definition~\ref{def:boundary-nondegeneracy} therefore imposes only the first-order normal expansion in \eqref{eq:boundary-local-expansion}, together with a positive definite quadratic term in the tangential variables.
\end{remark}

\begin{definition}[Simple regular GOE edge]\label{def:simple-regular-goe-edge}
Suppose that $g_\beta(\lambda_\beta)=0$. We say that $\lambda_\beta$ exposes a \textit{simple regular GOE edge} if the following properties hold.
\begin{enumerate}
\item The spectral distribution of $2\beta\mathsf M-D_{\lambda_\beta}\otimes1$ under $\tau_d$ has a density $\varrho_\beta$ near $0$, and $0$ is a regular right edge: $0\in\partial\operatorname{supp}\varrho_\beta$, there is $\delta>0$ such that $(0,\delta)\cap\operatorname{supp}\varrho_\beta=\varnothing$, and the square-root limit satisfies
\begin{equation*}
0<\lim_{\operatorname{supp}\varrho_\beta\ni x\uparrow0}
\frac{\varrho_\beta(x)}{\sqrt{|x|}}<\infty.
\end{equation*}
This is the regular-edge condition of \cite[Definition 2.4]{alt-erdos-kruger-schroder}, specialized to a right edge at $0$.
\item There are numbers $p_1,\ldots,p_d>0$ with $\sum_{a=1}^dp_a=1$ and $p_a:=-\partial_{\lambda_a}g_\beta(\lambda_\beta)$ for every $a\in\{1,\ldots,d\}$. With probability tending to one as $N\to\infty$, the largest eigenvalue of $2\beta M_N-D_{\lambda_\beta}\otimes I_N$ is simple. On this event, if $v_N$ is a unit eigenvector corresponding to this eigenvalue, then, for every fixed $\varepsilon>0$,
\begin{equation*}
\max_{a\in\{1,\ldots,d\}}
\left|\|P_av_N\|^2-p_a\right|
=O_{\mathbb P}(N^{-1/2+\varepsilon}).
\end{equation*}
\item There exists $\gamma_\beta>0$ such that
\begin{equation}\label{eq:simple-edge-tw}
\gamma_\beta N^{2/3}g_{N,\beta}(\lambda_\beta)
\xrightarrow[N\to\infty]{\mathrm d}\operatorname{TW}_1.
\end{equation}
\end{enumerate}
\end{definition}

Definition~\ref{def:simple-regular-goe-edge} records more random-matrix structure than is used directly in Theorem~\ref{thm:simple-regular-goe-free-energy}; the transfer theorem uses \eqref{eq:simple-edge-tw}, while the first two properties describe the deterministic edge and the limiting squared block norms $p_1,\ldots,p_d$. For correlated real symmetric ensembles satisfying the hypotheses of \cite[Theorem 2.11]{alt-erdos-kruger-schroder}, the regular-edge notion is \cite[Definition 2.4]{alt-erdos-kruger-schroder}, and GOE edge universality follows from \cite[Theorem 2.11]{alt-erdos-kruger-schroder}. Under the hypotheses of \cite[Theorem 2.6]{erdos-henheik-riabov}, its stochastic-domination estimate gives the second property at a regular edge when the deterministic microcanonical profile equals $(p_1,\ldots,p_d)$. Structured models with deterministic zero blocks can instead be checked directly.

\begin{lemma}[Boundary-value sensitivity]\label{lem:boundary-sensitivity}
Fix $\lambda_\beta\in\partial\Lambda_\beta$ satisfying Definition~\ref{def:boundary-multiplier} with $\theta_\beta>0$, and suppose that $(g_\beta,\widehat\phi_\beta)$ is nondegenerate at $\lambda_\beta$ in the sense of Definition~\ref{def:boundary-nondegeneracy}. Assume that $g_{N,\beta}(\lambda_\beta)=O_{\mathbb P}(N^{-2/3})$ and that
\begin{equation}\label{eq:single-edge-reduction}
F_N(\beta)
=\inf_{\substack{\lambda\in U\\g_\beta(\lambda)+g_{N,\beta}(\lambda_\beta)\leq0}}
\widehat\phi_\beta(\lambda)
+o_{\mathbb P}(N^{-2/3}).
\end{equation}
Then
\begin{equation}\label{eq:regular-edge-expansion}
F_N(\beta)-F(\beta)
=\theta_\beta g_{N,\beta}(\lambda_\beta)
+o_{\mathbb P}(N^{-2/3}).
\end{equation}
\end{lemma}

\begin{proof}
Use the parametrization in Definition~\ref{def:boundary-nondegeneracy} of the nondegeneracy of $(g_\beta,\widehat\phi_\beta)$ at $\lambda_\beta$. We take the neighborhood $U$ in \eqref{eq:single-edge-reduction} sufficiently small that its coordinates lie in a bounded neighborhood of $(0,0)$. Let $c_\beta>0$ satisfy $H_\beta\succeq c_\beta I_{d-1}$, and write
\begin{equation*}
\mathcal R_\beta(y,r):=
\widehat\phi_\beta(\Gamma(y,r))
-F(\beta)+\theta_\beta r-\frac12y^{\mathsf T}H_\beta y.
\end{equation*}
By \eqref{eq:boundary-local-expansion}, we have
\begin{equation}\label{eq:boundary-remainder}
\mathcal R_\beta(y,r)=o\bigl(|r|+\|y\|^2\bigr)
\end{equation}
as $(y,r)\to(0,0)$. Shrink the choice of $U$ before imposing \eqref{eq:single-edge-reduction}, if necessary, so that
\begin{equation}\label{eq:boundary-fixed-coercivity}
|\mathcal R_\beta(y,r)|
\leq\varepsilon_0\bigl(|r|+\|y\|^2\bigr),
\qquad
\varepsilon_0<\min\{\theta_\beta/2,c_\beta/4\},
\end{equation}
throughout $U$.

For deterministic $s$ sufficiently close to zero, the constraint $g_\beta(\lambda)+s\leq0$ becomes $r\leq-s$. Taking $y=0$ and $r=-s$, and then using \eqref{eq:boundary-remainder}, gives
\begin{equation}\label{eq:boundary-value-upper}
\inf_{\substack{\lambda\in U\\g_\beta(\lambda)+s\leq0}}
\widehat\phi_\beta(\lambda)
\leq F(\beta)+\theta_\beta s+o(|s|).
\end{equation}
For the reverse bound, choose $(y_s,r_s)$ satisfying $r_s\leq-s$ whose value is within $o(|s|)$ of the infimum. Equations \eqref{eq:boundary-value-upper} and \eqref{eq:boundary-fixed-coercivity} imply
\begin{equation}\label{eq:boundary-coordinate-localization}
|r_s|+\|y_s\|^2=O(|s|).
\end{equation}
Indeed, if $r_s\leq0$, then \eqref{eq:boundary-fixed-coercivity} gives
\begin{equation*}
\widehat\phi_\beta(\Gamma(y_s,r_s))-F(\beta)
\geq(\theta_\beta-\varepsilon_0)|r_s|
+\left(\frac{c_\beta}{2}-\varepsilon_0\right)\|y_s\|^2,
\end{equation*}
whereas if $r_s>0$, then $s<0$ and $0<r_s\leq|s|$; substituting $0<r_s\leq|s|$ into \eqref{eq:boundary-fixed-coercivity} gives $\|y_s\|^2=O(|s|)$. Combining \eqref{eq:boundary-coordinate-localization} with \eqref{eq:boundary-remainder} yields
\begin{align*}
\widehat\phi_\beta(\Gamma(y_s,r_s))-F(\beta)
=-\theta_\beta r_s+\frac12y_s^{\mathsf T}H_\beta y_s+o(|s|)\geq\theta_\beta s+o(|s|).
\end{align*}
Consequently,
\begin{equation}\label{eq:boundary-value-expansion}
\inf_{\substack{\lambda\in U\\g_\beta(\lambda)+s\leq0}}
\widehat\phi_\beta(\lambda)
=F(\beta)+\theta_\beta s+o(|s|).
\end{equation}

Since $g_{N,\beta}(\lambda_\beta)=O_{\mathbb P}(N^{-2/3})$, substituting $s=g_{N,\beta}(\lambda_\beta)$ in \eqref{eq:boundary-value-expansion} gives
\begin{equation*}
\inf_{\substack{\lambda\in U\\g_\beta(\lambda)+g_{N,\beta}(\lambda_\beta)\leq0}}
\widehat\phi_\beta(\lambda)
=F(\beta)+\theta_\beta g_{N,\beta}(\lambda_\beta)
+o_{\mathbb P}(N^{-2/3}).
\end{equation*}
Substitution into \eqref{eq:single-edge-reduction} proves \eqref{eq:regular-edge-expansion}.
\end{proof}

The edge law \eqref{eq:simple-edge-tw} is an assumption of the following transfer theorem; it is not deduced here from the finite fourth-moment assumption in \eqref{eq:coupling-covariance}.

\begin{theorem}[Free-energy fluctuations from one regular edge]\label{thm:simple-regular-goe-free-energy}
Fix $\lambda_\beta\in\partial\Lambda_\beta$ satisfying Definition~\ref{def:boundary-multiplier} with $\theta_\beta>0$. Suppose that $(g_\beta,\widehat\phi_\beta)$ is nondegenerate at $\lambda_\beta$ in the sense of Definition~\ref{def:boundary-nondegeneracy}, that $\lambda_\beta$ exposes a simple regular GOE edge in the sense of Definition~\ref{def:simple-regular-goe-edge}, and that the one-edge reduction \eqref{eq:single-edge-reduction} holds. If $\gamma_\beta$ is the edge-scaling constant in \eqref{eq:simple-edge-tw}, then
\begin{equation}\label{eq:simple-regular-edge-output}
\frac{\gamma_\beta}{\theta_\beta}N^{2/3}
\bigl(F_N(\beta)-F(\beta)\bigr)
\xrightarrow[N\to\infty]{\mathrm d}\operatorname{TW}_1.
\end{equation}
\end{theorem}

\begin{proof}
Definition~\ref{def:simple-regular-goe-edge} gives $g_{N,\beta}(\lambda_\beta)=O_{\mathbb P}(N^{-2/3})$ and \eqref{eq:simple-edge-tw}. Lemma~\ref{lem:boundary-sensitivity} therefore gives \eqref{eq:regular-edge-expansion}, and \eqref{eq:regular-edge-expansion} together with \eqref{eq:simple-edge-tw} proves \eqref{eq:simple-regular-edge-output}.
\end{proof}

\begin{remark}[The one-edge reduction and other boundary fluctuations]
The two-sided continuation in \eqref{eq:two-sided-objective} is essential in \eqref{eq:single-edge-reduction}. The random variable $N^{2/3}g_{N,\beta}(\lambda_\beta)$ has a nondegenerate limit with both signs. When $g_{N,\beta}(\lambda_\beta)<0$, the displaced constraint in \eqref{eq:single-edge-reduction} reaches points with $g_\beta(\lambda)>0$, which lie outside the closure of $\Lambda_\beta$.

Definition~\ref{def:simple-regular-goe-edge} controls $g_{N,\beta}(\lambda_\beta)$, the quantities $\|P_av_N\|^2$ for a unit eigenvector $v_N$ corresponding to the largest eigenvalue, and the limit in \eqref{eq:simple-edge-tw}; it does not imply the asymptotic identity \eqref{eq:single-edge-reduction}. To prove \eqref{eq:single-edge-reduction} for a general block pencil, one must additionally control the eigenvalues of $2\beta M_N-D_{\lambda_\beta}\otimes I_N$ within $O(N^{-2/3})$ of $0$, the corresponding quantities $\|P_av_N\|^2$, and the product-sphere integral as $\lambda_{\min}(D_\lambda\otimes I_N-2\beta M_N)$ approaches $0$.

If more than one eigenvalue of $2\beta M_N-D_{\lambda_\beta}\otimes I_N$ lies within $O_{\mathbb P}(N^{-2/3})$ of $0$, the $N^{2/3}$ fluctuation can depend on their joint limiting statistics. For a direct sum of independent scalar spherical models, it is a sum of independent Tracy--Widom fluctuations.
\end{remark}

\subsection{Ground-state fluctuations at a regular edge}

Set
\begin{equation*}
\mathsf{GS}_N:=\frac1N\max_{\sigma\in\mathcal S_N}H_N(\sigma).
\end{equation*}
Let $\mathsf{GS}$ denote the deterministic limit in \eqref{eq:ground-state-limit}. The objective in \eqref{eq:positive-spectral-dual} is invariant under $\lambda\mapsto\lambda+c\1$, so its minimization may be restricted to $\sum_{a=1}^d\lambda_a=0$. Let $\lambda_*$ be a minimizer on this hyperplane and define
\begin{equation*}
M_{*,N}:=M_N-D_{\lambda_*}\otimes I_N
\qquad\text{and}\qquad
e_*:=\lambda_{\max}(\mathsf M-D_{\lambda_*}\otimes1).
\end{equation*}

\begin{theorem}\label{thm:ground-state-fluctuations}
Suppose that the following two properties hold.
\begin{enumerate}
\item There exists $\gamma_*>0$ such that
\begin{equation}\label{eq:edge-tw-input}
\gamma_*N^{2/3}\bigl(\lambda_{\max}(M_{*,N})-e_*\bigr)
\xrightarrow[N\to\infty]{\mathrm d}\operatorname{TW}_1.
\end{equation}
\item With probability tending to one as $N\to\infty$, $M_{*,N}$ has a unit eigenvector $v_N$ corresponding to its largest eigenvalue and satisfying
\begin{equation}\label{eq:balanced-edge-vector}
\max_{a\in\{1,\ldots,d\}}
\left|\|P_av_N\|^2-\frac1d\right|
=o_{\mathbb P}(N^{-1/3}).
\end{equation}
\end{enumerate}
Then
\begin{equation}\label{eq:ground-state-tw}
\frac{\gamma_*}{d}N^{2/3}(\mathsf{GS}_N-\mathsf{GS})
\xrightarrow[N\to\infty]{\mathrm d}\operatorname{TW}_1.
\end{equation}
\end{theorem}

Property \eqref{eq:edge-tw-input} follows from \cite[Theorem 2.11]{alt-erdos-kruger-schroder} when $e_*$ is a regular right edge and the associated correlated real symmetric ensemble satisfies the assumptions of that theorem. Suppose that the hypotheses of \cite[Theorem 2.6]{erdos-henheik-riabov} hold and that its deterministic microcanonical block profile equals $(d^{-1},\ldots,d^{-1})$. At a regular edge, that theorem yields, for every fixed $\varepsilon>0$,
\begin{equation*}
\max_{a\in\{1,\ldots,d\}}
\left|\|P_av_N\|^2-\frac1d\right|
=O_{\mathbb P}(N^{-1/2+\varepsilon}),
\end{equation*}
so choosing $0<\varepsilon<1/6$ gives \eqref{eq:balanced-edge-vector}. In a differentiable setting with a simple exposed edge and a unit eigenvector $v$ corresponding to that edge, set $p_a:=\|P_av\|^2$. The first-order stationarity condition for the objective in \eqref{eq:positive-spectral-dual}, restricted to $\sum_a\lambda_a=0$, gives $1-dp_a=c$ for every $a\in\{1,\ldots,d\}$ and some $c\in\R$. Since $\sum_a p_a=1$, summing these identities gives $c=0$, and hence $p_a=d^{-1}$ for every $a\in\{1,\ldots,d\}$.

\begin{proof}
Write $e_N:=\lambda_{\max}(M_{*,N})$. On the event in the second assumption of the theorem, choose a unit eigenvector $v_N$ corresponding to $e_N$ and put $p_{N,a}:=\|P_av_N\|^2$. Equation \eqref{eq:balanced-edge-vector} implies that $p_{N,a}>0$ for every $a\in\{1,\ldots,d\}$ with probability tending to one as $N\to\infty$. On this event, define
\begin{equation*}
u_N:=\sum_{a=1}^d\frac{P_av_N}{\sqrt{p_{N,a}}}.
\end{equation*}
Then $\|P_au_N\|=1$ for every $a\in\{1,\ldots,d\}$. If $\widetilde v_N:=u_N-\sqrt d\,v_N$, then \eqref{eq:balanced-edge-vector} and the local Lipschitz continuity of the square root at $d^{-1}$ give
\begin{equation*}
\|\widetilde v_N\|^2
=\sum_{a=1}^d\left(1-\sqrt{dp_{N,a}}\right)^2
=o_{\mathbb P}(N^{-2/3}).
\end{equation*}
Proposition~\ref{prop:strong-convergence} gives $\|M_{*,N}\|=O_{\mathbb P}(1)$, and \eqref{eq:edge-tw-input} gives $e_N=O_{\mathbb P}(1)$. Therefore $\|M_{*,N}-e_NI_{dN}\|=O_{\mathbb P}(1)$. Since $M_{*,N}v_N=e_Nv_N$ and $\|u_N\|^2=\|\sqrt d\,v_N\|^2=d$, we obtain
\begin{equation}\label{eq:edge-energy-correction}
u_N^{\mathsf T}M_{*,N}u_N
=de_N+\widetilde v_N^{\mathsf T}(M_{*,N}-e_NI_{dN})\widetilde v_N
=de_N+o_{\mathbb P}(N^{-2/3}).
\end{equation}
Since $\|P_au_N\|=1$ for every $a\in\{1,\ldots,d\}$, the definition of $\mathsf{GS}_N$ and \eqref{eq:edge-energy-correction} give
\begin{equation}\label{eq:ground-state-lower-fluctuation}
\mathsf{GS}_N\geq\sum_{a=1}^d(\lambda_*)_a+de_N-o_{\mathbb P}(N^{-2/3}).
\end{equation}
Conversely, the inequality in Step~1 of the proof of Lemma~\ref{lem:spectral-dual}, applied with $M_N$ replaced by $M_{*,N}$, gives
\begin{equation}\label{eq:ground-state-upper-fluctuation}
\mathsf{GS}_N\leq\sum_{a=1}^d(\lambda_*)_a+de_N.
\end{equation}
By \eqref{eq:positive-spectral-dual},
\begin{equation*}
\mathsf{GS}=\sum_{a=1}^d(\lambda_*)_a+de_*.
\end{equation*}
Consequently, \eqref{eq:ground-state-lower-fluctuation} and \eqref{eq:ground-state-upper-fluctuation} imply
\begin{equation*}
\mathsf{GS}_N-\mathsf{GS}=d(e_N-e_*)+o_{\mathbb P}(N^{-2/3}).
\end{equation*}
Equation \eqref{eq:edge-tw-input} proves \eqref{eq:ground-state-tw}.
\end{proof}

The single-edge assumption in Theorem~\ref{thm:ground-state-fluctuations} is essential. If $H_N$ is a sum of independent scalar 2-spin spherical Hamiltonians, then $\mathsf{GS}_N$ is the corresponding sum of their ground-state energies. Its $N^{2/3}$ fluctuation is therefore a sum of independent Tracy--Widom variables rather than one Tracy--Widom variable.

\subsection{The bipartite spherical model}\label{sec:bipartite-fluctuations}

In this subsection, assume additionally that $g_{11}$ is Gaussian and consider the specialization in \eqref{eq:bipartite-specialization}. For $0<\beta<1$, \cite[Theorem 2.4(i)]{baik-lee} expresses the free energy through a logarithmic linear statistic evaluated a positive distance from the Marchenko--Pastur edge, with an error $O(N^{-2+\varepsilon})$. Thus Corollary~\ref{cor:gaussian-interior-fluctuations} gives the $N^{-1}$ Gaussian fluctuation scale, while \cite[Theorem 2.2(i)]{baik-lee} gives the explicit centering and variance. We now verify the hypotheses of Theorem~\ref{thm:simple-regular-goe-free-energy} for $\beta>1$. The critical value $\beta=1$ is excluded from both arguments.

For $d=2$ and $A_1=E_{12}$, we have
\begin{equation*}
M_N=\frac12
\begin{pmatrix}
0&G_{1,N}\\
G_{1,N}^{\mathsf T}&0
\end{pmatrix}.
\end{equation*}
Let $\mu_{\mathrm{MP}}$ be the equal-size Marchenko--Pastur law. The polar decomposition of $c_1$ and spectral calculus give
\begin{equation}\label{eq:bipartite-two-sided-potential}
\widehat\phi_\beta(\lambda)
=\frac{\lambda_1+\lambda_2-2}{2}
-\frac12\int_0^4\log\left|\lambda_1\lambda_2-\beta^2x\right|\,\mu_{\mathrm{MP}}(\mathrm dx),
\end{equation}
whenever the logarithm is integrable, and
\begin{equation}\label{eq:bipartite-edge-function}
g_\beta(\lambda)
=-\frac{\lambda_1+\lambda_2}{2}
+\sqrt{\left(\frac{\lambda_1-\lambda_2}{2}\right)^2+4\beta^2}.
\end{equation}
For $\lambda\in\Lambda_\beta$, one has $\lambda_1\lambda_2>4\beta^2$. For fixed $\lambda_1\lambda_2$, the first term in \eqref{eq:bipartite-two-sided-potential} is minimized at $\lambda_1=\lambda_2$. The Marchenko--Pastur Stieltjes transform
\begin{equation}\label{eq:mp-stieltjes}
\int_0^4\frac{\mu_{\mathrm{MP}}(\mathrm dx)}{z-x}
=\frac{1-\sqrt{1-4/z}}{2},
\qquad z\in\mathbb C\setminus[0,4],
\end{equation}
then gives, for $t>2\beta$ and $\beta>1$,
\begin{align*}
\frac{\mathrm d}{\mathrm dt}\widehat\phi_\beta(t,t)
=1-\frac{t}{\beta^2}\int_0^4\frac{\mu_{\mathrm{MP}}(\mathrm dx)}{t^2/\beta^2-x}
\geq1-\frac1\beta>0.
\end{align*}
Since $\log(4-x)$ is integrable against $\mu_{\mathrm{MP}}$, letting $t\downarrow2\beta$ in \eqref{eq:bipartite-two-sided-potential} yields
\begin{equation*}
\lambda_\beta=(2\beta,2\beta)\in\partial\Lambda_\beta
\qquad\text{and}\qquad
\widehat\phi_\beta(\lambda_\beta)=F(\beta).
\end{equation*}

Put $y:=(\lambda_1-\lambda_2)/2$ and $r:=g_\beta(\lambda)$. Using \eqref{eq:bipartite-edge-function} and the boundary value of \eqref{eq:mp-stieltjes} at $4$, one obtains, as $(y,r)\to(0,0)$,
\begin{equation}\label{eq:bipartite-boundary-expansion}
\widehat\phi_\beta(\lambda)
=F(\beta)-\frac{\beta-1}{\beta}r+\frac{y^2}{4\beta}+o\bigl(|r|+y^2\bigr).
\end{equation}
Equation \eqref{eq:bipartite-boundary-expansion} shows that $\theta_\beta$ is the boundary multiplier from Definition~\ref{def:boundary-multiplier} and that $(g_\beta,\widehat\phi_\beta)$ is nondegenerate at $\lambda_\beta$ in the sense of Definition~\ref{def:boundary-nondegeneracy}, with
\begin{equation*}
\theta_\beta=\frac{\beta-1}{\beta}>0
\qquad\text{and}\qquad
H_\beta=\frac1{2\beta}.
\end{equation*}

At $\lambda_\beta$, the shifted matrix is
\begin{equation*}
2\beta M_N-D_{\lambda_\beta}\otimes I_N
=
\begin{pmatrix}
-2\beta I_N&\beta G_{1,N}\\
\beta G_{1,N}^{\mathsf T}&-2\beta I_N
\end{pmatrix},
\end{equation*}
so
\begin{equation}\label{eq:bipartite-random-edge}
g_{N,\beta}(\lambda_\beta)=\beta\bigl(\|G_{1,N}\|-2\bigr).
\end{equation}
Choose unit vectors $u_N,v_N\in\R^N$ such that $G_{1,N}v_N=\|G_{1,N}\|u_N$ and $G_{1,N}^{\mathsf T}u_N=\|G_{1,N}\|v_N$; equivalently, $u_N$ and $v_N$ are left and right singular vectors associated with the largest singular value $\|G_{1,N}\|$. Then $(u_N,v_N)/\sqrt2$ is a unit eigenvector corresponding to the largest eigenvalue of $2\beta M_N-D_{\lambda_\beta}\otimes I_N$, and the squared norms of its two $N$-dimensional blocks are both $1/2=-\partial_{\lambda_a}g_\beta(\lambda_\beta)$. The Marchenko--Pastur density in \cite[Equation (2.31)]{baik-lee} vanishes as a positive constant times $\sqrt{4-x}$ as $x\uparrow4$, so $4$ is a regular right edge in the sense of Definition~\ref{def:simple-regular-goe-edge}(1). Moreover, \cite[Proposition 7.2]{baik-lee} gives
\begin{equation}\label{eq:bipartite-singular-tw}
2^{2/3}N^{2/3}\bigl(\|G_{1,N}\|-2\bigr)
\xrightarrow[N\to\infty]{\mathrm d}\operatorname{TW}_1.
\end{equation}
Hence $\lambda_\beta$ exposes a simple regular GOE edge in the sense of Definition~\ref{def:simple-regular-goe-edge}, with $\gamma_\beta=2^{2/3}/\beta$.

It remains to verify \eqref{eq:single-edge-reduction}. By \cite[Theorem 2.4(ii)]{baik-lee} and \eqref{eq:bipartite-specialization}, for any fixed $0<\varepsilon<1/100$,
\begin{equation}\label{eq:bipartite-free-edge-expansion}
F_N(\beta)-F(\beta)
=(\beta-1)\bigl(\|G_{1,N}\|-2\bigr)+o_{\mathbb P}(N^{-2/3}),
\end{equation}
where we used \eqref{eq:bipartite-singular-tw}. Equation \eqref{eq:boundary-value-expansion}, together with \eqref{eq:bipartite-random-edge} and $\theta_\beta=(\beta-1)/\beta$, gives the right-hand side of \eqref{eq:bipartite-free-edge-expansion} for the infimum in \eqref{eq:single-edge-reduction}. Thus \eqref{eq:single-edge-reduction} holds. Theorem~\ref{thm:simple-regular-goe-free-energy} now gives
\begin{equation*}
\frac{2^{2/3}}{\beta-1}N^{2/3}\bigl(F_N(\beta)-F(\beta)\bigr)
\xrightarrow[N\to\infty]{\mathrm d}\operatorname{TW}_1,
\end{equation*}
which agrees with \cite[Theorem 2.2(ii)]{baik-lee} after the rescaling in \eqref{eq:bipartite-specialization}.

For the ground state under \eqref{eq:bipartite-specialization}, we have $\mathsf{GS}_N=\|G_{1,N}\|$, so its fluctuations follow directly from \eqref{eq:bipartite-singular-tw}. This is also a special case of Theorem~\ref{thm:ground-state-fluctuations}.

\appendix
\section{Proof of Proposition~\ref{prop:strong-convergence}}\label{app:fourth-moment}

Recall $M_N$ from \eqref{eq:block-pencil} and $L_N(B,R,t)$ and $L(B,R,t)$ from \eqref{eq:affine-pencils}. We prove \eqref{eq:affine-spectrum} and \eqref{eq:affine-trace} under \eqref{eq:coupling-covariance}. For each fixed truncation level, we compare the independent matrix summands with Gaussian summands having the same entry covariances. We then remove truncation in operator norm, taking $N\to\infty$ before the truncation level tends to infinity.

\subsection{Comparison inequalities}

Let $H=H_0+\sum_\alpha Z_\alpha\in M_D(\C)$ be self-adjoint, where $H_0$ is deterministic and the $Z_\alpha$ are independent, centered, bounded, and self-adjoint. Define
\begin{align*}
\mathcal R_H&:=\max_\alpha\operatorname*{ess\,sup}\|Z_\alpha\|,
\quad
\sigma_H^2:=\left\|\sum_\alpha\E Z_\alpha^2\right\|,\quad 
\text{and}\quad\sigma_{*,H}^2&:=\sup_{\|u\|=\|v\|=1}\sum_\alpha\E|u^*Z_\alpha v|^2.
\end{align*}
Let $H^{\mathrm G}$ be a self-adjoint Gaussian matrix with the same mean and the same joint covariance of the real and imaginary entries as $H$. Brailovskaya and van Handel \cite[Theorem 2.6]{brailovskaya-van-handel} give, for every $s\geq0$,
\begin{gather}
\P\left(d_{\mathrm H}(\operatorname{spec}H,\operatorname{spec}H^{\mathrm G})>C\varepsilon_H(s)\right)\leq De^{-s},\label{eq:comparison-spectrum}\\
\varepsilon_H(s):=\sigma_{*,H}\sqrt{s}+\mathcal R_H^{1/3}\sigma_H^{2/3}s^{2/3}+\mathcal R_Hs.\label{eq:comparison-error}
\end{gather}
Here $C$ is universal, and \eqref{eq:comparison-spectrum} holds for every joint realization of $H$ and $H^{\mathrm G}$; see the paragraph following that theorem. In particular, these matrices may be independent.

For $p\in\mathbb N$, write $a_p(H):=(D^{-1}\E\Tr H^{2p})^{1/(2p)}$. Brailovskaya and van Handel \cite[Theorem 2.9 with $q=\infty$, and Lemma 9.20]{brailovskaya-van-handel} also give
\begin{equation}\label{eq:comparison-moments}
|a_p(H)-a_p(H^{\mathrm G})|\leq C_p\bigl(\mathcal R_H^{1/3}\sigma_H^{2/3}+\mathcal R_H\bigr),
\end{equation}
and, for every $s\geq p$,
\begin{equation}\label{eq:comparison-moment-concentration}
\P\left(\left|(D^{-1}\Tr H^{2p})^{1/(2p)}-a_p(H)\right|
>C_p\left[(\sigma_{*,H}+\mathcal R_H^\frac{1}{2}a_p(H)^\frac{1}{2})\sqrt{s}+\mathcal R_Hs\right]\right)
\leq2e^{-s}.
\end{equation}
The hypotheses of these inequalities require independence of the summands $Z_\alpha$, with no independence assumption on the entries within each summand.

\subsection{Bounded couplings}

For $K\in\mathbb N$, define
\begin{equation}\label{eq:fourth-truncation}
g_{ij}^{(K)}:=g_{ij}\mathbbm{1}_{\{\|g_{ij}\|_{\mathrm F}\leq K\}}
-\E\left[g_{11}\mathbbm{1}_{\{\|g_{11}\|_{\mathrm F}\leq K\}}\right],
\end{equation}
and let $M_N^{(K)}$ be \eqref{eq:block-pencil} with $g_{ij}$ replaced by $g_{ij}^{(K)}$. Thus $\E g_{ij}^{(K)}=0$, $\|g_{ij}^{(K)}\|_{\mathrm F}\leq2K$, and the matrices $g_{ij}^{(K)}$ remain independent over ordered pairs $(i,j)$.

\begin{proof}[Proof of \eqref{eq:affine-spectrum}--\eqref{eq:affine-trace} for $M_N^{(K)}$]
\noindent\emph{Step 1: Bounds for $\mathcal R_H$, $\sigma_H$, and $\sigma_{*,H}$.}
Fix $B=B^*\in M_d(\C)$, $R\in M_d(\C)$, and $t\in\R$. For $i\leq j$, set
\begin{equation*}
U_{ij}^{(K)}:=\frac12\bigl(g_{ij}^{(K)}+(g_{ji}^{(K)})^{\mathsf T}\bigr)
\qquad\text{and}\qquad
\widetilde U_{ij}^{(K)}:=RU_{ij}^{(K)}R^*.
\end{equation*}
For $i<j$, define
\begin{equation*}
Z_{ij}:=\frac{t}{\sqrt N}\left(\widetilde U_{ij}^{(K)}\otimes E_{ij}
+(\widetilde U_{ij}^{(K)})^*\otimes E_{ji}\right),
\qquad\text{and}\qquad
Z_{ii}:=\frac{t}{\sqrt N}\widetilde U_{ii}^{(K)}\otimes E_{ii}.
\end{equation*}
Then
\begin{equation}\label{eq:bounded-pencil-sum}
L_N^{(K)}(B,R,t):=B\otimes I_N+t(R\otimes I_N)M_N^{(K)}(R^*\otimes I_N)
=B\otimes I_N+\sum_{i\leq j}Z_{ij}.
\end{equation}
The $Z_{ij}$ are independent and centered self-adjoint matrices. In particular, for $i<j$, the entire pair $(g_{ij}^{(K)},g_{ji}^{(K)})$ occurs in one summand. The diagonal variables $g_{ii}^{(K)}$ are independent of these pairs.

Put $c_{K,R}:=\max\{\E\|\widetilde U_{12}^{(K)}\|^2,\E\|\widetilde U_{11}^{(K)}\|^2\}$. For $H:=L_N^{(K)}(B,R,t)$ in \eqref{eq:bounded-pencil-sum}, we have
\begin{equation}\label{eq:bounded-comparison-parameters}
\mathcal R_H\leq\frac{2|t|\|R\|^2K}{\sqrt N},\qquad
\sigma_H^2\leq t^2c_{K,R},
\qquad\text{and}\qquad
\sigma_{*,H}^2\leq\frac{4t^2c_{K,R}}N.
\end{equation}
Indeed, $\|Z_{ij}\|=|t|N^{-1/2}\|\widetilde U_{ij}^{(K)}\|$. For $i<j$, we have
\begin{equation*}
Z_{ij}^2=\frac{t^2}N\left(
\widetilde U_{ij}^{(K)}(\widetilde U_{ij}^{(K)})^*\otimes E_{ii}
+(\widetilde U_{ij}^{(K)})^*\widetilde U_{ij}^{(K)}\otimes E_{jj}\right).
\end{equation*}
Together with the corresponding diagonal term, this gives $\|\sum_{i\leq j}\E Z_{ij}^2\|\leq t^2c_{K,R}$. Finally, write $u=(u_i)_{i=1}^N$ and $v=(v_i)_{i=1}^N$ in site coordinates, with $u_i,v_i\in\C^d$. For $i<j$, we have
\begin{equation*}
u^*Z_{ij}v=\frac{t}{\sqrt N}\left(u_i^*\widetilde U_{ij}^{(K)}v_j
+u_j^*(\widetilde U_{ij}^{(K)})^*v_i\right).
\end{equation*}
Using $|x+y|^2\leq2|x|^2+2|y|^2$ and summing the off-diagonal and diagonal terms gives
\begin{equation*}
\sum_{i\leq j}\E|u^*Z_{ij}v|^2
\leq\frac{2t^2c_{K,R}}N\sum_{i\neq j}\|u_i\|^2\|v_j\|^2
+\frac{t^2c_{K,R}}N\sum_i\|u_i\|^2\|v_i\|^2
\leq\frac{4t^2c_{K,R}}N
\end{equation*}
when $\|u\|=\|v\|=1$. This proves \eqref{eq:bounded-comparison-parameters}.

\smallskip
\noindent\emph{Step 2: Gaussian matrices with the same covariance.}
Let $\Delta^{(K)}$ be the covariance matrix of the vector $(g_{11}^{(K),ab})_{a,b}$, and let $A_r^{(K)}\in M_d(\R)$ be obtained by reshaping column $r$ of $(\Delta^{(K)})^\frac{1}{2}$, for $r\in\{1,\ldots,d^2\}$. Thus
\begin{equation*}
\Delta_{ab,ce}^{(K)}=\sum_{r=1}^{d^2}(A_r^{(K)})_{ab}(A_r^{(K)})_{ce}.
\end{equation*}
On a product extension of the probability space, take an infinite array $(\widetilde z_{ij}^r)_{i,j\in\mathbb N,\,r\in\{1,\ldots,d^2\}}$ of independent standard real Gaussian variables, independent of $(g_{ij})_{i,j\in\mathbb N}$, and use this one array for every $N$ and $K$. The Gaussian matrices $\sum_{r=1}^{d^2}A_r^{(K)}\widetilde z_{ij}^r$ have the same covariance as $g_{ij}^{(K)}$ and are independent over ordered pairs. Writing
\begin{equation*}
\widetilde G_{r,N}:=N^{-1/2}(\widetilde z_{ij}^r)_{i,j\in\{1,\ldots,N\}},
\end{equation*}
define directly
\begin{multline*}
L_N^{\mathrm G,(K)}(B,R,t):=B\otimes I_N
+\frac t2\sum_{r=1}^{d^2}\bigl[(RA_r^{(K)}R^*)\otimes\widetilde G_{r,N}
+(R(A_r^{(K)})^{\mathsf T}R^*)\otimes\widetilde G_{r,N}^{\mathsf T}\bigr].
\end{multline*}
This matrix has the same full entry covariance as $L_N^{(K)}(B,R,t)$, including the diagonal site blocks, and is therefore an admissible $H^{\mathrm G}$ in \eqref{eq:comparison-spectrum} and \eqref{eq:comparison-moments} for $H=L_N^{(K)}(B,R,t)$.

Let $\widetilde c_1,\ldots,\widetilde c_{d^2}$ be a fixed free standard circular family, and put
\begin{equation}\label{eq:truncated-free-matrix}
\mathsf M^{(K)}:=\frac12\sum_{r=1}^{d^2}\left(A_r^{(K)}\otimes\widetilde c_r+(A_r^{(K)})^{\mathsf T}\otimes\widetilde c_r^*\right).
\end{equation}
Xiang and Zhang \cite[Theorem 3.3, Remark 3.6, and Corollary 3.7]{xiang-zhang}, applied to the Gaussian array just constructed, give \eqref{eq:affine-spectrum} and \eqref{eq:affine-trace} for $L_N^{\mathrm G,(K)}$ with limit
\begin{equation*}
L^{(K)}(B,R,t):=B\otimes1+t(R\otimes1)\mathsf M^{(K)}(R^*\otimes1).
\end{equation*}
The Gaussian norm bound \eqref{eq:ginibre-norm-moments} also gives, for every fixed $q<\infty$,
\begin{equation}\label{eq:gaussian-comparison-integrability}
\sup_N\E\|L_N^{\mathrm G,(K)}(B,R,t)\|^q<\infty.
\end{equation}
Thus every fixed normalized trace moment of the Gaussian matrices converges both almost surely and in expectation.

\smallskip
\noindent\emph{Step 3: Convergence for fixed $K$.}
Apply \eqref{eq:comparison-spectrum} with $D=dN$ and $s=4\log(dN)$. Equations \eqref{eq:comparison-error} and \eqref{eq:bounded-comparison-parameters} give
\begin{equation*}
\varepsilon_{L_N^{(K)}}(4\log(dN))
=O_{K,R,t}\left(N^{-1/2}\sqrt{\log N}+N^{-1/6}(\log N)^{2/3}+N^{-1/2}\log N\right).
\end{equation*}
This tends to zero as $N\to\infty$. Since $\sum_N(dN)^{-3}<\infty$, \eqref{eq:comparison-spectrum}, Borel--Cantelli, and Step 2 imply
\begin{equation}\label{eq:bounded-spectrum-limit}
\lim_{N\to\infty}d_{\mathrm H}\bigl(\operatorname{spec}L_N^{(K)}(B,R,t),\operatorname{spec}L^{(K)}(B,R,t)\bigr)=0,
\end{equation}
almost surely on the product space. The event in \eqref{eq:bounded-spectrum-limit} depends only on the original array $(g_{ij})$ and a deterministic limit. Fubini's theorem therefore gives probability one for this event on the original probability space.

Fix $p\in\mathbb N$. Equations \eqref{eq:comparison-moments} and \eqref{eq:gaussian-comparison-integrability} give
\begin{equation*}
\lim_{N\to\infty}a_p(L_N^{(K)}(B,R,t))
=\bigl(\tau_d(L^{(K)}(B,R,t)^{2p})\bigr)^{1/(2p)}.
\end{equation*}
The threshold in \eqref{eq:comparison-moment-concentration}, with $s=4\log(dN)$, is $O_{K,B,R,t,p}(N^{-1/4}\sqrt{\log N}+N^{-1/2}\log N)$. Its exceptional probabilities are at most $2(dN)^{-4}$ and are summable. Borel--Cantelli therefore gives
\begin{equation}\label{eq:bounded-even-moment-limit}
\lim_{N\to\infty}\frac1{dN}\Tr L_N^{(K)}(B,R,t)^{2p}
=\tau_d(L^{(K)}(B,R,t)^{2p}).
\end{equation}
Apply \eqref{eq:bounded-even-moment-limit} with $B$ replaced by $B+sI_d$ for $2p+1$ distinct rational values of $s$. Expanding $(L_N^{(K)}+sI_{dN})^{2p}$ yields an invertible Vandermonde system for its moments of orders $0,\ldots,2p$. Hence every normalized trace moment converges to the corresponding moment of $L^{(K)}$.

Intersect these probability-one events over $K,p\in\mathbb N$ and over $B=B^*$, $R$, and $t$ with rational real and imaginary coordinates. For each fixed $K$, \eqref{eq:bounded-spectrum-limit} with $(B,R,t)=(0,I_d,1)$ gives an eventual bound on $\|M_N^{(K)}\|$. The equicontinuity argument in Step~1 of the proof of Lemma~\ref{lem:spectral-dual} then extends the spectral and polynomial-trace limits locally uniformly to all coefficients. All spectra lie in one compact interval on each compact coefficient set. Polynomial approximation on that interval gives, on this one event, for every continuous $f:\R\to\C$,
\begin{equation}\label{eq:bounded-trace-limit}
\lim_{N\to\infty}\frac1{dN}\Tr f(L_N^{(K)}(B,R,t))
=\tau_d(f(L^{(K)}(B,R,t))),
\end{equation}
locally uniformly in $(B,R,t)$ for each fixed $f$.
\end{proof}

\subsection{Removal of truncation}

Recall $g_{ij}^{(K)}$ from \eqref{eq:fourth-truncation} and set $\rho_K:=(\E\|g_{11}-g_{11}^{(K)}\|_{\mathrm F}^2)^\frac{1}{2}$. Centering in $L^2$ gives
\begin{equation}\label{eq:fourth-tail-variance}
\rho_K^2\leq\E\left[\|g_{11}\|_{\mathrm F}^2\mathbbm{1}_{\{\|g_{11}\|_{\mathrm F}>K\}}\right],
\qquad\text{and}\qquad
\lim_{K\to\infty}\rho_K=0.
\end{equation}
For each fixed $K,a,b$, the entries $(g_{ij}-g_{ij}^{(K)})_{ab}$ form an infinite i.i.d. centered real array with finite fourth moment. For
\begin{equation*}
Y_{ab,N}^{(K)}:=N^{-1/2}\bigl((g_{ij}-g_{ij}^{(K)})_{ab}\bigr)_{i,j\in\{1,\ldots,N\}},
\end{equation*}
Bai and Yin \cite[Theorem 2.1, with $k=1$]{bai-yin-products} give
\begin{equation}\label{eq:coordinate-tail-norm}
\limsup_{N\to\infty}\|Y_{ab,N}^{(K)}\|
\leq2\bigl(\E|(g_{11}-g_{11}^{(K)})_{ab}|^2\bigr)^\frac{1}{2}\leq2\rho_K,
\end{equation}
almost surely. If the variance is zero, the corresponding array vanishes almost surely, so \eqref{eq:coordinate-tail-norm} also holds. There are only finitely many coordinates and countably many $K$, so these bounds hold simultaneously. Symmetrizing as in \eqref{eq:block-pencil} yields
\begin{equation}\label{eq:fourth-tail-norm}
\begin{aligned}
\limsup_{N\to\infty}\|M_N-M_N^{(K)}\|
&\leq\sum_{a,b=1}^d\limsup_{N\to\infty}\|Y_{ab,N}^{(K)}\|
\stackrel{\eqref{eq:coordinate-tail-norm}}{\leq}2d^2\rho_K,\\
\lim_{K\to\infty}\limsup_{N\to\infty}\|M_N-M_N^{(K)}\|&=0.
\end{aligned}
\end{equation}
No independence between different coordinates $(a,b)$ is used in this argument.

Let $A_r^{(\infty)}$ be obtained by reshaping column $r$ of $\Delta^\frac{1}{2}$, and define $\mathsf M^{(\infty)}$ by \eqref{eq:truncated-free-matrix} with these coefficients and the same circular variables. Equation \eqref{eq:fourth-tail-variance} gives $\lim_{K\to\infty}\|\Delta^{(K)}-\Delta\|=0$. Continuity of the positive semidefinite square root therefore implies
\begin{equation}\label{eq:free-truncation-limit}
\lim_{K\to\infty}\|\mathsf M^{(K)}-\mathsf M^{(\infty)}\|
\leq\lim_{K\to\infty}\sum_{r=1}^{d^2}\|A_r^{(K)}-A_r^{(\infty)}\|\,\|\widetilde c_r\|=0.
\end{equation}
The operator $\mathsf M^{(\infty)}$ has the same $M_d(\C)$-valued distribution as $\mathsf M$, since both are semicircular with covariance map \eqref{eq:eta}. Faithfulness of the traces identifies the spectra of their self-adjoint polynomials: the spectrum is the support of the scalar spectral measure. Thus $\mathsf M^{(\infty)}$ may be used in place of $\mathsf M$ in \eqref{eq:affine-spectrum}--\eqref{eq:affine-trace}.

For fixed $B,R,t$, the triangle inequality and the norm bound for perturbations of self-adjoint spectra give
\begin{align*}
&\limsup_{N\to\infty}d_{\mathrm H}\bigl(\operatorname{spec}L_N(B,R,t),\operatorname{spec}L(B,R,t)\bigr)\leq |t|\|R\|^2\left(2d^2\rho_K+\|\mathsf M^{(K)}-\mathsf M^{(\infty)}\|\right),
\end{align*}
where \eqref{eq:bounded-spectrum-limit} removes the middle term. Letting $K\to\infty$ and using \eqref{eq:fourth-tail-variance} and \eqref{eq:free-truncation-limit} proves \eqref{eq:affine-spectrum}.

For clarity, all trace comparisons can be made on a common compact spectral interval. Fixing one $K$ in \eqref{eq:fourth-tail-norm} first gives $\limsup_N\|M_N\|<\infty$, and
\begin{equation*}
\limsup_{N\to\infty}\|M_N^{(K)}\|
\leq\limsup_{N\to\infty}\|M_N\|+2d^2\rho_K.
\end{equation*}
Together with \eqref{eq:free-truncation-limit}, this bounds the spectra of all four matrices involved in the comparison, uniformly for sufficiently large $K$ and then sufficiently large $N$. On a compact coefficient set $\mathcal C$, choose an interval $[-T,T]$ containing these spectra, and let $p$ be a polynomial. For self-adjoint $A,B$ with norms at most $T$, telescoping powers gives
\begin{equation*}
\|p(A)-p(B)\|\leq C_{p,T}\|A-B\|.
\end{equation*}
The trace bound $|D^{-1}\Tr X|\leq\|X\|$, its analogue for $\tau_d$, and \eqref{eq:bounded-trace-limit} imply
\begin{multline*}
\limsup_{N\to\infty}\sup_{(B,R,t)\in\mathcal C}
\left|\frac1{dN}\Tr f(L_N(B,R,t))-\tau_d(f(L(B,R,t)))\right|\\
\leq4\|f-p\|_{\infty,[-T,T]}
+C_{p,T}\sup_{(B,R,t)\in\mathcal C}|t|\|R\|^2
\left(2d^2\rho_K+\|\mathsf M^{(K)}-\mathsf M^{(\infty)}\|\right).
\end{multline*}
First let $K\to\infty$, and then approximate $f$ uniformly by polynomials on $[-T,T]$. This proves \eqref{eq:affine-trace} locally uniformly in the coefficients. The countable intersection in the bounded proof and \eqref{eq:fourth-tail-norm} give one probability-one event, independent of $B,R,t,f$. Proposition~\ref{prop:strong-convergence} follows.

\bibliographystyle{amsplain}
\bibliography{ref}

\end{document}